\documentclass[12pt]{amsart}
\usepackage{amssymb,amsmath,amscd}

\calclayout

\newtheorem{theorem}{Theorem}
\newtheorem{corollary}{Corollary}
\newtheorem{proposition}{Proposition}
\newtheorem{lemma}{Lemma}

\newtheorem{defn}{Definition}

\newcommand{\cF}{{\mathcal F}}
\newcommand{\cG}{{\mathcal G}}
\newcommand{\cI}{{\mathcal I}}

\newcommand{\cO}{{\mathcal O}}
\newcommand{\cX}{{\mathcal X}}

\newcommand{\mA}{{\mathbb A}}
\newcommand{\mG}{{\mathbb G}}
\newcommand{\mP}{{\mathbb P}}
\newcommand{\mZ}{{\mathbb Z}}

\newcommand{\diag}{\operatorname{diag}}
\newcommand{\di}{\operatorname{div}}
\newcommand{\gr}{\operatorname{gr}}
\newcommand{\Spec}{\operatorname{Spec}}

\title{A geometric approach to Standard Monomial Theory}
\author{M.~Brion}
\address{Institut Fourier\\
UMR 5582 du CNRS\\
F-38402 Saint-Martin d'H\`eres Cedex}
\email{Michel.Brion@ujf-grenoble.fr}
\author{V.~Lakshmibai}
\thanks{Second author partially supported by N.S.F. Grant DMS-9971295}
\address{Department of Mathematics\\
Northeastern University\\
Boston, MA 02115-5096\\}
\email{lakshmibai@neu.edu} 

\date{}
\begin{document}

\begin{abstract}
We obtain a geometric construction of a ``standard monomial basis''
for the homogeneous coordinate ring associated with any ample line
bundle on any flag variety. This basis is compatible with Schubert
varieties, opposite Schubert varieties, and unions of intersections of
these varieties. Our approach relies on vanishing theorems and a
degeneration of the diagonal ; it also yields a standard monomial basis
for the multi--homogeneous coordinate rings of flag varieties of
classical type.
\end{abstract}

\maketitle

\section*{Introduction}

Consider the Grassmannian $X$ of linear subspaces of dimension $r$ in
$k^n$, where $k$ is a field. We regard $X$ as a closed subvariety of
projective space $\mP(\wedge^r k^n)$ via the Pl\"ucker embedding; let
$L$ be the corresponding very ample line bundle on $X$. Then the
ring $\bigoplus_{m=0}^{\infty} H^0(X,L^{\otimes m})$
admits a nice basis, defined as follows.

\smallskip

Let $\{v_1,\ldots,v_n\}$ be the usual basis of $k^n$; then the
$v_{i_1}\wedge\cdots\wedge v_{i_r}$, $1\le i_1 < \cdots < i_r\le n$,
form a basis of $\wedge^r k^n$. We put $I=(i_1,\ldots, i_r)$,
$v_I = v_{i_1} \wedge \cdots \wedge v_{i_r}$, and we denote by
$\{p_I\}$ the dual basis of the basis $\{v_I\}$; the $p_I$ 
(regarded in $H^0(X,L)$) are the Pl\"ucker coordinates. Define a
partial order on the set $\cI$ of indices $I$ by letting
$I=(i_1,\ldots,i_r)\le (j_1,\ldots,j_r)=J$ if and only if 
$i_1\le j_1,\ldots,i_r\le j_r$. Then 

\smallskip

\noindent
(i) {\it The monomials $p_{I_1}p_{I_2}\cdots p_{I_m}$ where
$I_1,\ldots,I_m\in\cI$ satisfy $I_1\le I_2\le \cdots\le I_m$, form a
basis of $H^0(X,L^{\otimes m})$.} 

\smallskip

\noindent
(ii) {\it For any $I,J\in\cI$, we have
$p_I\, p_J- \sum_{I',J',\;I'\le I,J\le J'} \, a_{I'J'}\, p_{I'}\, p_{J'}
=0,$
where $a_{I'J'}\in k$.}

\smallskip

The monomials in (i) are called the {\it standard monomials of degree}
$m$, and the relations in (ii) are the {\it quadratic straightening
relations}; they allow to express any non--standard monomial in the
$p_I$ as a linear combination of standard monomials.

\smallskip

Further, this {\it standard monomial basis} of the homogeneous
coordinate ring of $X$ is compatible with its Schubert subvarieties,
in the following sense. For any $I\in\cI$, let 
$X_I=\{V\in X ~\vert~
\dim(V\cap\text{span}(v_1,\ldots,v_s)) \ge \#(j,i_j\le s),
1\le s\le r\}$
be the corresponding Schubert variety; then the restriction
$p_J\vert_{X_I}$ is nonzero if and only if $J\le I$. The monomial
$p_{I_1}\cdots p_{I_m}$ will be called {\it standard on} $X_I$ if
$I_1\le\cdots\le I_m\le I$; equivalently, this monomial is standard
and does not vanish identically on $X_I$. Now 

\smallskip

\noindent
(iii) {\it The standard monomials of degree $m$ on $X_I$ restrict to a
basis of $H^0(X_I,L^{\otimes m})$. The standard monomials of degree
$m$ that are not standard on $X_I$, form a basis of the kernel of the
restriction map 
$H^0(X,L^{\otimes m})\rightarrow H^0(X_I,L^{\otimes m})$.}

\smallskip

These classical results go back to Hodge, see \cite{H}. They have
important geometric consequences, e.g., $X$ is projectively normal in
the Pl\"ucker embedding; its homogeneous ideal is generated by the
quadratic straightening relations; the homogeneous ideal of any
Schubert variety $X_I$ is generated by these relations together with
the $p_J$ where $J\not\le I$. 

\smallskip

The purpose of {\it Standard Monomial Theory} (SMT) is to generalize
Hodge's results to any flag variety $X=G/P$ (where $G$ is a semisimple
algebraic group over an algebraically closed field $k$, and $P$ a
parabolic subgroup) and to any effective line bundle $L$ on $X$. 
SMT was developed by Lakshmibai, Musili, and Seshadri in a series of
papers, culminating in \cite{LS1} where it is established for all
classical groups $G$. There the approach goes by ascending induction
on the Schubert varieties, using their partial resolutions as
projective line bundles over smaller Schubert varieties.

\smallskip

Further results concerning certain exceptional or Kac--Moody groups
led to conjectural formulations of a general SMT, see \cite{LS2}.
These conjectures were then proved by Littelmann, who introduced
new combinatorial and algebraic tools: the path model of
representations of any Kac--Moody group, and Lusztig's Frobenius map
for quantum groups at roots of unity (see \cite{Li1,Li2}). 

\smallskip

In the present paper, we obtain a geometric construction of
a SMT basis for $H^0(X,L)$, where $X=G/P$ is any flag variety
and $L$ is any ample line bundle on $X$. This basis is compatible
with Schubert varieties (that is, with orbit closures in $X$ of a
Borel subgroup $B$ of $G$) and also with opposite Schubert
varieties (the orbit closures of an opposite Borel subgroup
$B^-$); in fact, it is compatible with any intersection of a
Schubert variety with an opposite Schubert variety. We call such
intersections {\it Richardson varieties}, since they were first
considered by Richardson in \cite{Ric}. Our approach adapts to the
case where $L$ is an effective line bundle on a flag variety
{\it of classical type} in the sense of \cite{LS1}. This sharpens
the results of \cite{LS1} concerning the classical groups.

\smallskip

Our work may be regarded as one step towards a purely geometric
proof of Littelmann's results concerning SMT. He constructed a
basis of $T$--eigenvectors for $H^0(X,L)$ (where $T$ is the maximal
torus common to $B$ and $B^-$) indexed by certain piecewise linear
paths in the positive Weyl chamber, called {\it LS paths}. This basis 
turns out to be compatible with Richardson varieties; notice that
these are $T$--invariant. In fact, the endpoints of the path indexing a
basis vector parametrize the smallest Richardson variety where this
vector does not vanish identically (see \cite{LL}). If $L$ is
associated with a weight of classical type, then the LS paths are just
line segments: they are uniquely determined by their endpoints. This
explains a posteriori why our geometric approach completes the program 
of SMT in that case.

\smallskip

In fact, our approach of SMT for an ample line bundle
$L$ on a flag variety $X$ uses little of the rich geometry
and combinatorics attached to $X$. Specifically, we only rely
on vanishing theorems for unions of Richardson varieties
(these being direct consequences of the existence of a
Frobenius splitting of $X$, compatible with Schubert varieties
and opposite Schubert varieties), together with the following
property.

\smallskip

\noindent
(iv) {\it The diagonal in $X\times X$ admits a flat $T$--invariant
degeneration to the union of all products $X_w\times X^w$,
where the $X_w$ are the Schubert varieties and the $X^w$ are
the corresponding opposite Schubert varieties.}

\smallskip

The latter result follows from \cite{BP} (we provide a direct
proof in Section 3). It plays an essential r\^ole in establishing
generalizations of (i) and (iii); conversely, it turns out that 
the existence of a SMT basis implies (iv), see the Remark after
Proposition 7.

\smallskip

It is worth noticing that (iv) is a stronger form
of the fact that the classes of Schubert varieties form a
free basis of the homology group (or Chow group) of $X$,
the dual basis for the intersection pairing consisting of
the classes of opposite Schubert varieties. This fact
(in a different formulation) has been used by Knutson to
establish an asymptotic version of the Littelmann character
formula, see \cite{Kn}.

\medskip

This paper is organized as follows. In the preliminary Section
1, we introduce notation and study the geometry of Richardson
varieties. Vanishing theorems for cohomology groups
of line bundles on Richardson varieties are established in Section
2, by slight generalizations of the methods of Frobenius splitting.
In Section 3, we construct filtrations of the $T$--module
$H^0(X,L)$ that are compatible with restrictions to Richardson
varieties. Our SMT basis of $H^0(X,L)$ is defined in Section 4;
it is shown to be compatible with all unions of Richardson
varieties. In Section 5, we generalize statements (i) and (iii)
above to any ample line bundle $L$ on a flag variety $G/P$; then
(ii) follows from (i) together with compatibility properties of
our basis. The case where the homogeneous line bundle $L$ is
associated with a weight of classical type
(e.g., a fundamental weight of a classical group)
is considered in detail in Section 6. There we give a geometric
characterization of the {\it admissible pairs} of \cite{LS1} 
(these parametrize the weights of the $T$--module $H^0(X,L)$).
The final Section 7 develops SMT for those effective
line bundles that correspond to sums of weights of classical type.

\section{Richardson varieties}

The ground field $k$ is algebraically closed, of arbitrary
characteristic. Let $G$ be a simply--connected semisimple algebraic
group. Choose opposite Borel subgroups $B$ and $B^-$ of $G$, with
common torus $T$; let $\cX(T)$ be the group of characters of $T$, also
called weights. In the root system $R$ of $(G,T)$, we have the subset
$R^+$ of positive roots (that is, of roots of $(B,T)$), and the subset
$S$ of simple roots. For each $\alpha\in R$, let $\check\alpha$ be the
corresponding coroot and let $U_\alpha$ be the corresponding additive
one--parameter subgroup of $G$, normalized by $T$.

We also have the Weyl group $W$ of $(G,T)$; for each $\alpha\in R$, we
denote by $s_\alpha\in W$ the corresponding reflection. Then the
group $W$ is generated by the simple reflections $s_\alpha$,
$\alpha\in S$; this defines the length function $\ell$ and the Bruhat
order $\leq$ on $W$. Let $w_o$ be the longest element of $W$, then
$B^-=w_o B w_o$.

Let $P$ be a parabolic subgroup of $G$ containing $B$ and let $W_P$ be
the Weyl group of $(P,T)$, a parabolic subgroup of $W$; let $w_{o,P}$
be the longest element of $W_P$. Each right $W_P$--coset in $W$
contains a unique element of minimal length; this defines the subset
$W^P$ of minimal representatives of the quotient $W/W_P$. This subset
is invariant under the map $w\longmapsto w_oww_{o,P}$; the induced
bijection of $W^P$ reverses the Bruhat order. 

Each character $\lambda$ of $P$ defines a $G$--linearized line bundle on
the homogeneous space $G/P$; we denote that line bundle by $L_\lambda$.
The assignement $\lambda\longmapsto L_\lambda$ yields an isomorphism from
the character group $\cX(P)$ to the Picard group of $G/P$. Further,
the line bundle $L_\lambda$ is generated by its global sections if and
only if $\lambda$ (regarded as a character of $T$) is dominant; in
that case, $H^0(G/P,L_\lambda)$ is a $G$--module with lowest weight
$-\lambda$. 

Let $W_\lambda$ be the isotropy group of $\lambda$ in $W$, and let
$P_\lambda$ be the parabolic subgroup of $G$ generated by $B$ and 
$W_\lambda$; then $W_\lambda\supseteq W_P$, $W^\lambda\subseteq W^P$,
and $P_\lambda\supseteq P$. We shall identify $W^\lambda$ with the
$W$--orbit of the weight $\lambda$, and denote by $w(\lambda)$ the
image of $w\in W$ in $W/W_\lambda\simeq W^\lambda$.

The {\it extremal weight vectors} 
$p_{w(\lambda)}\in H^0(G/P,L_\lambda)$ are the $T$--eigenvectors of
weight $-w(\lambda)$ for some $w\in W^\lambda$. These vectors are
uniquely defined up to scalars.

We say that $\lambda$ is $P$--{\it regular} if $P_\lambda=P$. The ample
line bundles on $G/P$ are the $L_\lambda$ where $\lambda$ is dominant
and $P$--regular; under these assumptions, $L_\lambda$ is in fact very
ample. We may then identify each $w\in W^P$ to $w(\lambda)$, and we put
$p_w=p_{w(\lambda)}$.

The $T$--fixed points in $G/P$ are the $e_w=wP/P$
($w\in W/W_P$); we index them by $W^P$. The $B$--orbit $C_w=Be_w$ is a
{\it Bruhat cell}, an affine space of dimension $\ell(w)$; it
closure in $G/P$ is the {\it Schubert variety} $X_w$. The complement 
$X_w - C_w$ is the {\it boundary} $\partial X_w$. We have
$$
\partial X_w=\bigcup_{v\in W^P,\;v<w} X_v,
$$
and the irreducible components of $\partial X_w$ are the 
{\it Schubert divisors} $X_v$ where $v\in W^P$, $v<w$ and
$\ell(v)=\ell(w)-1$. Then there exists $\beta\in R^+$ such that
$v=ws_\beta$.

Let $\lambda$ be a character of $P$ and let $f_w$ be the restriction
to $X_w$ of the natural map $G/P\longrightarrow G/P_\lambda$; then
$f_w(X_w)=X_{w(\lambda)}$. The set 
$$
\partial_\lambda X_w:=f_w^{-1}(\partial X_{w(\lambda)})
$$
is called the $\lambda$--{\it boundary} of $X_w$; it is the union of the
Schubert divisors $X_{ws_\beta}$ where 
$\langle\lambda,\check\beta\rangle > 0$. If
$\lambda$ is dominant, then we have by Chevalley's formula:
$$
\di(p_{w(\lambda)}\vert_{X_w})=
\sum \langle\lambda,\check\beta\rangle X_{ws_\beta}
$$
(sum over all $\beta\in R^+$ such that $X_{s_\beta w}$ is a divisor
in $X_w$). In particular, the zero set of $p_{w(\lambda)}$ in $X_w$ is
$\partial_\lambda X_w$. If in addition $\lambda$ is $P$--regular, then
$\partial_\lambda X_w = \partial X_w$.

We shall also need the {\it opposite Bruhat cell} $C^w=B^- e_w$ of
codimension $\ell(w)$ in $G/P$, the {\it opposite Schubert variety}
$X^w$ (the closure of $C^w$) and its boundary $\partial X^w$. Then 
$X^w = w_oX_{w_oww_{o,P}}$ and 
$$
\partial X^w=\bigcup_{v\in W^P, v>w} X^v.
$$

Recall that all Schubert varieties are normal and Cohen--Macaulay
(thus, the same holds for all opposite Schubert varieties). Further,
all scheme--theoretic intersections of unions of Schubert varieties and
opposite Schubert varieties are reduced (see \cite{RR,Ram1,Ram2}).

\begin{defn}
Let $v,w$ in $W^P$. We call the intersection
$$
X_w^v:=X_w\cap X^v
$$ 
a Richardson variety in $G/P$. We define its boundaries by
$$
(\partial X_w)^v:=\partial X_w\cap X^v \text{ and } 
(\partial X^v)_w:=X_w\cap \partial X^v.
$$
\end{defn}

Notice that $X_w^v$ and its boundaries are closed reduced, $T$--stable
subschemes of $G/P$. The $X_w^v$ were considered by Richardson,
who showed e.g. that they are irreducible (see \cite{Ric}; the
intersections $C_w\cap C^v$ were analyzed by Deodhar, see
\cite{D}). We shall give another proof of this result, and obtain a
little more.

\begin{lemma}\label{richardson}

\begin{enumerate}

\item $X_w^v$ is non--empty if and only if $v\leq w$; then $X_w^v$ is
irreducible of dimension $\ell(w)-\ell(v)$, and $(\partial X_w)^v$,
$(\partial X^v)_w$ have pure codimension $1$ in $X_w^v$. Further,
$X_w^v$ is normal and Cohen--Macaulay. 

\item The $T$--fixed points in $X_w^v$ are the $e_x$ where $x\in W^P$
and $v\leq x\leq w$.

\item For $x,y$ in $W^P$, we have 
$X_y^x\subseteq X_w^v \Longleftrightarrow v\leq x\leq y\leq w$.
\end{enumerate}

\end{lemma}

\begin{proof}
(2) is evident; it implies (3) and the first assertion of (1). 
To prove the remaining assertions, we use a variant of the argument of
\cite{B} Lemma 2. Consider the fiber product $G\times^B X_w$ with
projection map
$$
p:G\times^B X_w\longrightarrow G/B,
$$
a $G$--equivariant locally trivial fibration with fiber $X_w$. We also
have the ``multiplication'' map
$$
m:G\times^B X_w\longrightarrow G/P,~(g,x)\longmapsto gx.
$$
This is a $G$--equivariant map to $G/P$; thus, it is also a locally
trivial fibration. Its fiber $m^{-1}(e_1)$ is isomorphic to
$\overline{Pw^{-1}B}/B$ (a Schubert variety in $G/B$).

Next let $i:X^v\longrightarrow G/P$ be the inclusion and consider the
cartesian product
$$
Z=X^v\times_{G/P} (G\times^B X_w)
$$
with projections $\iota$ to $G\times^B X_w$, $\mu$ to $X^v$ and $\pi$
to $G/B$, as
displayed in the following commutative diagram:
$$
\CD
 G/B @<{\pi}<< Z @>{\mu}>> X^v\\
 @V{id}VV    @V{\iota}VV @V{i}VV\\
 G/B @<{p}<<   G\times^B X_w @>{m}>> G/P
\endCD
$$
By definition, the square on the right is cartesian, so that $\mu$ is
also a locally trivial fibration with fiber $\overline{Pw^{-1}B}/B$
and base $X^v$. Since Schubert varieties are irreducible, normal and
Cohen--Macaulay, it follows that the same holds for $Z$. Further, we
have
$$
\dim(Z)=\dim(G\times^B X_w)+\dim(X^v)-\dim(G/P)
=\dim(G/B)+\ell(w)-\ell(v).
$$

Notice that the fiber of $\pi:Z\longrightarrow G/B$ at each $gB/B$
identifies to the intersection $X^v\cap gX_w$; in particular,
$\pi^{-1}(B/B)=X_w^v$. Notice also that $\iota:Z\longrightarrow
G\times^B X_w$ is a closed immersion with $B^-$--stable image (since
this holds for $i:X^v\longrightarrow G/P$). Thus, $B^-$ acts on $Z$ so
that $\pi$ is equivariant. Since $B^-B/B$ is an open neighborhood of
$B/B$ in $G/B$, isomorphic to $U^-$, its pullback under $\pi$ is an
open subset of $Z$, isomorphic to $U^-\times X_w^v$. Therefore,
$X_w^v$ is irreducible, normal and Cohen--Macaulay of dimension
$\ell(w)-\ell(v)$. 
\end{proof}

We also record the following easy result, to be used in Section 7.

\begin{lemma}\label{extremal}
Let $v\le w$ in $W^P$, let $\lambda$ be a dominant character of
$P$ and let $x(\lambda)\in W^\lambda$. Then the restriction of 
$p_{x(\lambda)}$ to $X_w^v$ is non--zero if and only if $x(\lambda)$
admits a lift $x\in W^P$ such that $v\le x\le w$. Further, the ring
$$
\bigoplus_{n=0}^{\infty} H^0(X_w^v,L_{n\lambda})
$$
is integral over its subring generated by the
$p_{x(\lambda)}\vert_{X_w^v}$ where $x\in W^P$ and $v\le x\le w$.
\end{lemma}

\begin{proof}
Consider the natural map $G/P\to G/P_\lambda$ and its restriction
$f:X_w^v\longrightarrow f(X_w^v)$. The open subset 
$(p_{x(\lambda)}\ne 0)$ of $G/P_\lambda$ is affine, $T$--stable and
contains $e_{x(\lambda)}$ as its unique closed $T$--orbit. Thus,
$p_{x(\lambda)}\vert_{X_w^v}\ne 0$ if and only if 
$e_{x(\lambda)}\in f(X_w^v)$. By Borel's fixed point theorem, this
amounts to the existence of a $T$--fixed point $e_x\in X_w^v$ such that
$f(e_x)=e_{x(\lambda)}$. Now Lemma \ref{richardson} (2) completes the
proof of the first assertion.

By the preceding arguments, the sections
$p_{x(\lambda)}\vert_{X_w^v}$, $x\in W^P$, $v\le x\le w$ do not vanish
simultaneously at a $T$--fixed point of $X_w^v$. Since these sections
are eigenvectors of $T$, it follows that they have no common
zeroes. This implies the second assertion.
\end{proof}

\noindent
{\bf Remark.} The image of a Richardson variety $X_w^v$ under a
morphism $G/P\longrightarrow G/P_\lambda$ need not be another
Richardson variety. Consider for example $G=SL(3)$ with
simple reflections $s_1$, $s_2$. Let $P=B$, $w=s_2 s_1$, $v=s_2$ and
$\lambda=\omega_1$ (the fundamental weight fixed by $s_2$). Then
$X_w^v$ is one--dimensional and mapped isomorphically to its image
$f(X_w^v)$ in $G/P_\lambda$. Since the $T$--fixed points in $f(X_w^v)$
are $e_{\omega_1}$ and $e_{s_2s_1(\omega_1)}$, it follows that
$f(X_w^v)$ is not a Richardson variety.

\section{Cohomology vanishing for Richardson varieties}

In this section, we assume that the characteristic of $k$ is
$p>0$. Let $X$ be a scheme of finite type over $k$. Let
$F:X\longrightarrow X$ be the absolute Frobenius morphism, that is,
$F$ is the identity map on the topological space of $X$, and
$F^{\#}:\cO_X\longrightarrow F_*\cO_X$ is the $p$--th power map. Then
$X$ is called {\it Frobenius split} if the map 
$F^{\#}$ is split. We shall need a slight generalization of this
notion, involving the composition $F^r = F\circ\cdots\circ F$ ($r$
times), where $r$ is any positive integer.

\begin{defn}
We say that $X$ is split if there exists a positive integer $r$ such
that the map 
$$
(F^r)^{\#}:\cO_X\longrightarrow F^r_*\cO_X
$$ 
splits, that is, there exists an $\cO_X$--linear map 
$$
\varphi:F^r_*\cO_X\longrightarrow\cO_X
$$ 
such that $\varphi\circ (F^r)^{\#}$ is the identity; then $\varphi$ is
called a splitting.
\end{defn}

We shall also need a slight generalization of the notion of 
Frobenius splitting relative to an effective Cartier divisor 
(see \cite{Ram2}).

\begin{defn}
Let $X$ be a normal variety and $D$ an effective Weil divisor on $X$,
with canonical section $s$. We say that $X$ is $D$--split if there
exist a positive integer $r$ and an $\cO_X$--linear map
$$
\psi:F^r_*\cO_X(D)\longrightarrow \cO_X
$$ 
such that the map 
$$
\varphi:F^r_*\cO_X\longrightarrow \cO_X,~f\longmapsto \psi(fs)
$$
is a splitting. Then $\psi$ is called a $D$--splitting.

We say that a closed subscheme $Y$ of $X$, with ideal sheaf $\cI_Y$,
is compatibly $D$--split if 
(a) no irreducible component of $Y$ is contained in the support of
$D$, and 
(b) $\varphi(F^r_*\cI_Y)=\cI_Y$.
\end{defn}

\noindent
{\bf Remarks.} (i) Let $U$ be an open subset of $X$ such that $X-U$
has codimension at least $2$ in $X$. Then $X$ is $D$--split if and only
if $U$ is $D\cap U$--split (to see this, let $i:U\longrightarrow X$ be the
inclusion, then $i_*\cO_U=\cO_X$ and $i_*\cO_U(D\cap U) = \cO_X(D)$ by
normality of $X$.)

Let $Y$ be a closed subscheme of $X$ such that $Y\cap U$ is dense in
$Y$. Then $Y$ is compatibly $D$--split if and only if $Y\cap U$ is
compatibly $D\cap U$--split (this is checked by the arguments of
\cite{Ram2} 1.4--1.7).

(ii) If $X$ is split compatibly with an effective Weil divisor $D$,
then $X$ is $(p^r-1)D$--split (to see this, one may assume that $X$ is
nonsingular, by (i). Let $\varphi$ be a compatible splitting, then
$\varphi(F^r_*\cO_X(-D))=\cO_X(-D)$. Define 
$\psi:F^r_*\cO_X(D)\longrightarrow \cO_X$ by
$\psi(f\sigma^{p^r-1})=\sigma\varphi(f\sigma^{-1})$ for any local
sections $f$ of $\cO_X$ and $\sigma$ of $\cO_X(D)$. Then one checks
that $\psi$ is well--defined, $\cO_X$--linear and satisfies 
$\psi(f s^{p^r-1})=\varphi(f)$. 

(iii) Let $D$ and $E$ be effective Weil divisors in $X$, such that
$D-E$ is effective. If $X$ is $D$--split, then it is $E$--split as well;
if in addition a closed subscheme $Y$ of $X$ is compatibly $D$--split,
then it is compatibly $E$--split 
(this follows from (i) together with \cite{Ram2} Remark 1.3 (ii).)

\begin{lemma}\label{iterate} 
Let $D$, $E$ be effective Weil divisors on a normal variety $X$, such
that the support of $D$ contains the support of $E$. If $X$ is
$D$--split, then $X$ is $E$--split as well. If moreover a closed
subscheme $Y$ of $X$ is compatibly $D$--split, then $X$ is compatibly
$E$--split.
\end{lemma}

\begin{proof}
Let $U$ be the set of those points of $X$ at which $D$ is a Cartier
divisor. Then $U$ is an open subset with complement of codimension at
least $2$ (since $U$ contains the nonsingular locus of $X$). Moreover,
$Y\cap U$ is dense in $Y$ (since $U$ contains the complement of the
support of $D$). Thus, by Remark (i), we may replace $X$ with $U$, and
hence assume that $D$ is a Cartier divisor.

Now let $\psi:F^r_*\cO_X(D)\longrightarrow \cO_X$ be a
$D$--splitting. We regard $\psi$ as an additive map
$\cO_X(D)\longrightarrow\cO_X$ such that: $\psi(s)=1$, and
$\psi(f^{p^r}\sigma)=f\psi(\sigma)$ for any local sections $f$ of
$\cO_X$ and $\sigma$ of $\cO_X(D)$. For any positive integer $n$, we
set
$$
{\bf n} = p^{r(n-1)}+p^{r(n-2)}+\cdots+1
$$
(then ${\bf 1}=1$), and we define inductively a map 
$$
\psi^n:F^{rn\#}\cO_X({\bf n}D)\longrightarrow\cO_X
$$
by: $\psi^1=\psi$, and 
$$
\psi^n(f\sigma^{{\bf n}})=\psi(\psi^{n-1}(f\sigma^{{\bf n-1}})\sigma)
$$
for any local sections $f$ of $\cO_X$ and $\sigma$ of $\cO_X(D)$. 
Then one may check that $\psi^n$ is well defined and is a 
${\bf n}D$--splitting of $X$. If moreover a closed subscheme $Y$ is
compatibly $D$--split, then $\psi(F^r_*(\cI_Y s))=\cI_Y$. By induction,
it follows that $\psi^n(F^{rn}_*(\cI_Y s^{{\bf n}}))=\cI_Y$, so that
$Y$ is compatibly ${\bf n}D$--split.

Since the support of $D$ contains the support of $E$, there exists a
positive integer $n$ such that ${\bf n}D - E$ is effective. Then $X$
is ${\bf n}D$--split, so that it is $E$-split by Remark (ii).
\end{proof}

\begin{lemma}\label{surjective}
Let $X$ be a normal projective variety endowed with an effective Weil
divisor $D$ and with a globally generated line bundle $L$; let $Y$ be
a closed subscheme of $X$. Assume that (a) $X$ is $D$--split compatibly
with $Y$, and (b) the support of $D$ contains the support of an
effective ample divisor. Then $H^i(X,L)=0=H^i(Y,L)$ for all $i\geq 1$,
and the restriction map $H^0(X,L)\longrightarrow H^0(Y,L)$ is surjective.
\end{lemma}

\begin{proof}
Choose an effective ample Cartier divisor $E$, with support contained
in the support of $D$. Then $X$ is $E$--split compatibly with $Y$, by
Lemma \ref{iterate}. Now the assertions follow from \cite{Ram2} 1.12,
1.13.
\end{proof}

We now apply this to Richardson varieties. By \cite{Ram2} 3.5, the
variety $G/P$ is split compatibly with all Schubert varieties and with
all opposite Schubert varieties; as a consequence, $G/P$ is split
compatibly with all unions of Richardson varieties. By \cite{Ram2}
1.10, it follows that all scheme--theoretical intersections of unions
of Richardson varieties are reduced; and using \cite{Ram2} 1.13, this
also implies

\begin{lemma}\label{vanishing}
Let $\lambda$ be a regular dominant character of $P$ and let $Z$ be a
union of Richardson varieties in $G/P$. Then the restriction map
$H^0(G/P,L_\lambda)\longrightarrow H^0(Z,L_\lambda)$ is surjective, and 
$H^i(Z,L_\lambda)=0$ for all $i\ge 1$. As a consequence,
$H^i(X,L_\lambda\otimes\cI_Z)=0$ for all $i\geq 1$.
\end{lemma}

\noindent
{\bf Remark.}
If we only assume that $\lambda$ is dominant, then Lemma \ref{vanishing}
extends to all unions of Schubert varieties (by \cite{Ram2}), but not to
all unions of Richardson varieties. As a trivial example, take
$G/P=\mP^1$, the projective line with $T$--fixed points 
$0$ and $\infty$, and $\lambda=0$. Then $Z:=\{0,\infty\}$ is a union
of Richardson varieties, and the restriction map 
$H^0(\mP^1,\cO_{\mP^1})\longrightarrow H^0(Z,\cO_Z)$ is not
surjective. As a less trivial example, take $G/P=\mP^1\times\mP^1$, 
$Z=(\mP^1\times\{0,\infty\})\cup(\{0,\infty\}\times\mP^1)$, and
$\lambda=0$. Then $Z$ is again a union of Richardson varieties, and
one checks that $H^1(Z,\cO_Z)\neq 0$.

However, Lemma \ref{vanishing} does extend to all dominant
characters and to unions of Richardson varieties 
{\it with a common index}.

\begin{proposition}\label{vanishing1}
Let $\lambda$ be a dominant character of $P$ and let $Z$ be a union of
Richardson varieties $X_w^v$ in $G/P$, all having the same $w$. Then
the restriction
$H^0(G/P,L_{\lambda})\rightarrow H^0(Z,L_{\lambda})$ 
is surjective, and $H^i(Z,L_{\lambda})=0$ for all $i\geq 1$.

As a consequence, we have
$H^i(X_w^v,L_\lambda(-Z))=0$ for all $i\geq 1$, where 
$v\leq w$ in $W^P$, and $Z$ is a union of irreducible components of
$(\partial X^v)_w$.
\end{proposition}

\begin{proof}
The Schubert variety $X_w$ is split compatibly with the effective Weil
divisor $\partial X_w$ and with $Z$. By assumption, $\partial X_w$
contains no irreducible component of $Z$. Using Remarks (i) and (ii),
it follows that $X_w$ is $(p-1)\partial X_w$--split compatibly with
$Z$. Further, $\partial X_w$ is the support of an ample effective
divisor, as follows from Chevalley's formula. Thus, Lemma
\ref{surjective} applies and yields surjectivity of 
$H^0(X_w,L_\lambda)\longrightarrow H^0(Z,L_\lambda)$ together with
vanishing of $H^i(Z,L_\lambda)$ for $i\geq 1$. Now surjectivity of
$H^0(G/P,L_\lambda)\longrightarrow H^0(X_w,L_\lambda)$ completes the
proof of the first assertion.

In particular, we have 
$H^i(X_w^v,L_\lambda)=H^i(Z,L_\lambda)=0$
for all $i\geq 1$, and the restriction map
$H^0(X_w^v,L_\lambda)\longrightarrow H^0(Z,L_\lambda)$
is surjective; this implies the second assertion.
\end{proof}

We shall also need the following, more technical vanishing result.

\begin{proposition}\label{vanishing2}
Let $\lambda$ be a dominant character of $P$ and let $v,w$ in
$W^P$ such that $v\leq w$. Then 
$$
H^i(X_w^v,L_\lambda(-(\partial_{\lambda} X_w)^v - Z))=0
$$
for any $i\geq 1$ and for any (possibly empty) union $Z$ of
irreducible components of $(\partial X^v)_w$.
\end{proposition}

\begin{proof}
We shall rely on the following result (see \cite{MK} Theorem 1). Let
$\pi:X\longrightarrow Y$ be a proper morphism of schemes. Let $D$
(resp.~$E$) be a closed subscheme of $X$ (resp.~$Y$) and let $i$ be a
positive integer such that:

\noindent
(i) $\pi^{-1}(E)$ is contained in $D$ (as sets).

\noindent
(ii) $R^i\pi_*(\cI_D)=0$ outside $E$.

\noindent
(iii) $X$ is split compatibly with $D$.

\noindent
Then $R^i\pi_*(\cI_D)=0$ everywhere.

To apply this result, consider the restriction 
$$
f:X_w^v\longrightarrow f(X_w^v)
$$
of the natural map $G/P\longrightarrow G/P_\lambda$. Then
$L_\lambda=f^*M_\lambda$ for a very ample line bundle $M_\lambda$ on
$f(X_w^v)$. Let $Y$ be the corresponding affine cone over $f(X_w^v)$,
with vertex $0$ and projection map
$$
q:Y-\{0\}\longrightarrow f(X_w^v).
$$ 
And let $X$ be the total space of the line bundle $L_{-\lambda}$ (dual
to $L_{\lambda}$), with projection map 
$$
p:X\longrightarrow X_w^v
$$ 
and zero section $X_0$. Then the algebra
$$
H^0(X,\cO_X)=\bigoplus_{n=0}^{\infty}\, H^0(X_w^v,L_{n\,\lambda})
$$
contains $H^0(Y,\cO_Y)$ as the subalgebra generated by
$H^0(f(X_w^v),M_\lambda)$. The algebra $H^0(X,\cO_X)$
is finitely generated, and the corresponding morphism
$$
X\longrightarrow \Spec\, H^0(X,\cO_X)
$$ 
is proper, since the line bundle $L_{\lambda}$ is globally generated. 
Moreover, since $L_{\lambda}$ is the pullback under $f$ of the
very ample line bundle $M_{\lambda}$, the algebra $H^0(X,\cO_X)$ is a
finite module over its subalgebra $H^0(Y,\cO_Y)$. This defines a
proper morphism $\pi:X\longrightarrow Y$, and we have
$\pi^{-1}(0)=X_0$ (as sets). Moreover, the diagram
$$
\CD
 X-X_0  @>{\pi}>> Y-\{0\}\\
 @V{p}VV    @V{q}VV\\
 X^v_w  @>{f}>> f(X_w^v)
\endCD
$$
is cartesian, and the vertical maps are principal $\mG_m$--bundles. 

Now let $D=X_0\cup p^{-1}((\partial_{\lambda} X_w)^v\cup Z)$; this is
a closed subscheme of $X$ with ideal sheaf 
$$
\cI_D=p^*L_{\lambda}(-(\partial_{\lambda} X_w)^v-Z).
$$
Let $E$ be the affine cone over $f((\partial_\lambda X_w)^v)$;
this is a closed subscheme of $Y$. We check that the conditions (i),
(ii) and (iii) hold.

For (i), notice that
$$
\pi^{-1}(E)=X_0\cup p^{-1}f^{-1}f((\partial_\lambda X_w)^v)
= X_0\cup p^{-1}((\partial_{\lambda} X_w)^v)
$$
(as sets), by the definition of $(\partial_{\lambda} X_w)^v$. In other
words, $\pi^{-1}(E)\subseteq D$ as sets.

For (ii), observe that $\cI_D=p^*L_\lambda(-Z)$ outside
$\pi^{-1}(E)$. Thus, (ii) is equivalent to: 
$R^i\pi_*(p^*L_\lambda(-Z))=0$ outside $E$. We show that 
$R^i\pi_*(p^*L_\lambda(-Z))=0$ outside $0$. Using the cartesian
square above, it suffices to check that $R^if_*(L_\lambda(-Z))=0$;
by the Leray spectral sequence and the Serre vanishing theorem, this
amounts to $H^i(X_w^v,L_{n\lambda}(-Z))=0$ for large $n$. But this
holds by Proposition \ref{vanishing1}.

For (iii), recall that $X_w^v$ is split compatibly with 
$(\partial_{\lambda} X_w)^v\cup Z$. Let $\varphi$ be a compatible
splitting; then $\varphi$ lifts uniquely to a splitting of $X$
compatibly with $X_0$ and with 
$p^{-1}((\partial_{\lambda} X_w)^v\cup Z)$.
It follows that $X$ is split compatibly with $D$.

We thus obtain: $R^i\pi_*(\cI_D)=0$ for all $i\geq 1$. Since $Y$ is
affine, this amounts to: $H^i(X,\cI_D)=0$ for all $i\geq 1$. On the
other hand, since the morphism $p:X\longrightarrow X_w^v$ is affine, 
we obtain that $H^i(X,\cI_D)$ is isomorphic to 
$$
H^i(X_w^v,p_*p^*(L_{\lambda}(-(\partial_{\lambda}X_w)^v)-Z))=
H^i(X_w^v,L_{\lambda}(-(\partial_{\lambda}X_w)^v-Z)\otimes p_*\cO_X).
$$
Further, $p_*\cO_X$ contains $\cO_{X_w^v}$ as a direct factor. This
yields
$$
H^i(X_w^v,L_{\lambda}(-(\partial_{\lambda}X_w)^v-Z))=0
$$
for all $i\geq 1$.
\end{proof}

\begin{corollary}\label{surj}
With the above notations, the restriction map
$$
H^0(X_w^v,L_{\lambda}(-(\partial_\lambda X_w)^v))\longrightarrow 
H^0(X_w^x,L_{\lambda}(-(\partial_\lambda X_w)^x))
$$
is surjective for any $x\in W^P$ such that $v\leq x\leq w$.
\end{corollary}

\begin{proof}
We may reduce to the case that $\ell(x)=\ell(v)+1$, that is, $X_w^x$
is an irreducible component of $(\partial X^v)_w$. Then 
$H^1(X_w^v,L_{\lambda}(-(\partial_\lambda X_w)^v-X_w^x))=0$
by Proposition \ref{vanishing2}. Now the assertion follows from the
exact sequence 
$$
0\longrightarrow 
L_{\lambda}\vert_{X_w^v}(-(\partial_\lambda X_w)^v-X_w^x)
\longrightarrow 
L_{\lambda}\vert_{X_w^v}(-(\partial_\lambda X_w)^v)
\longrightarrow 
L_{\lambda}\vert_{X_w^x}(-(\partial_\lambda X_w)^x)
\longrightarrow 0.
$$
\end{proof}

Notice finally that Lemma \ref{vanishing}, Propositions
\ref{vanishing1} and \ref{vanishing2}, and Corollary \ref{surj} also
hold in characteristic zero, as follows from the argument in
\cite{Ram2} 3.7. 

\section{Filtrations}

In this section, we shall obtain natural filtrations of the
$T$--modules $H^0(X_w^v,L_\lambda)$ and 
$H^0(X_w^v,L_\lambda(-(\partial_\lambda X_w)^v))$ (where $X_w^v$ is a
Richardson variety in $G/P$, and $\lambda$ is a dominant character of
$P$), and we shall describe their associated graded modules. For this,
we shall construct a degeneration of $X_w^v$ embedded diagonally in
$G/P\times G/P$, to a union of products of Richardson varieties.

Such a degeneration was obtained in \cite{BP} Theorem 16 for 
$X_w^v = G/B$, by using the wonderful compactification of the adjoint
group of $G$; it was extended to certain subvarieties in $G/P$,
including Schubert varieties, in \cite{B} Theorem 2. Here we follow a
direct, self--contained approach, at the cost of repeating some of the
arguments in \cite{BP} and \cite{B}. We begin by establishing a
K\"unneth decomposition of the class of the diagonal of $G/P$, in the 
Grothendieck group of $G/P\times G/P$; such a decomposition is deduced
in \cite{BP} from a degeneration of the diagonal.

Let $K(G/P\times G/P)$ be the Grothendieck group of the category of
coherent sheaves on $G/P\times G/P$. The class of a coherent sheaf
$\cF$ in this group will be denoted by $[\cF]$. 

\begin{lemma}\label{formula}
We have in $K(G/P\times G/P)$:
$$\displaylines{
[\cO_{\diag(G/P)}] = [\cO_{\bigcup_{x\in W^P} X_x\times X^x}]
\hfill\cr\hfill
=\sum_{x\in W^P} [\cO_{X_x}(-\partial X_x)\otimes\cO_{X^x}]
=\sum_{x\in W^P} [\cO_{X_x}\otimes\cO_{X^x}(-\partial X^x)].
\cr}$$
\end{lemma}

\begin{proof}
Let $Z=\bigcup_{x\in W^P} X_x\times X^x$.  We first claim that
$$
[\cO_Z]=\sum_{x\in W^P} [\cO_{X_x}(-\partial X_x)\otimes\cO_{X^x}].
$$
Let $W^P =\{x_1,\ldots,x_N\}$ be an indexing such that $i\le j$
whenever $x_i\le x_j$. Then one obtains easily:
$$
(X_{x_i}\times X^{x_i})\cap(\bigcup_{j<i} X_{x_j}\times X^{x_j})
=\partial X_{x_i}\times X^{x_i}.
$$
Now let $\cO_{Z,\ge i}$ be the subsheaf of $\cO_Z$ consisting of those
sections that vanish on $X_{x_j}\times X^{x_j}$ for each $j<i$. Then
the $\cO_{Z,\ge i}$ are a decreasing filtration of $\cO_Z$, and 
$$
\cO_{Z,\ge i}/\cO_{Z,\ge i+1}\simeq 
\cO_{Z,\ge i}\vert_{X_{x_i}\times X^{x_i}}
\simeq \cO_{X_i}(-\partial X_i)\otimes\cO_{X^i}.
$$
Further, $[\cO_Z]=\sum_{i=1}^N [\cO_{Z,\ge i}/\cO_{Z,\ge i+1}$]
in $K(G/P\times G/P)$. This implies our claim.

One checks similarly that
$$
[\cO_Z]=\sum_{x\in W^P} [\cO_{X_x}\otimes\cO_{X^x}(-\partial X^x)],
$$
using the increasing filtration of $\cO_Z$ by the subsheaves
$\cO_{Z,\le i}$ consisting of those sections that vanish on
$X_{x_j}\times X^{x_j}$ for each $j>i$.
 
To complete the proof, it suffices to check that 
$$
[\cO_{\diag(G/P)}] = 
\sum_{x\in W^P} [\cO_{X_x}\otimes\cO_{X^x}(-\partial X^x)].\eqno(*)
$$
For this, we recall some well--known facts on Grothendieck groups
of flag varieties.

Since the Bruhat cells $C_x$, $x\in W^P$, form a cellular 
decomposition of $G/P$,
the abelian group $K(G/P)$ is generated by the $[\cO_{X_x}]$, 
$x\in W^P$. Likewise, it is generated by the $[\cO_{X^y}]$, 
$y\in W^P$. Further, $K(G/P)$ is a ring for the product
$$
[\cF]\cdot[\cG]=\sum_{i\ge 0} (-1)^i[Tor_i^{G/P}(\cF,\cG)],
$$
and the Euler characteristic of coherent sheaves yields an additive map
$$
\begin{matrix}
\chi: & K(G/P) & \longrightarrow & \mZ \\
& [\cF] & \longmapsto & \chi(\cF). \\
\end{matrix}
$$
Since $X_x$ and $X^y$ are Cohen--Macaulay and intersect properly in
$G/P$, we have $Tor_i^{G/P}(\cO_{X_x},\cO_{X^y})=0$ for all 
$i\ge 1$ (see \cite{B} Lemma 1 for details). 
And since the intersection $X_x\cap X^y=X_x^y$ is reduced, we obtain 
$$
Tor_0^{G/P}(\cO_{X_x},\cO_{X^y})=
\cO_{X_x}\otimes_{\cO_{G/P}}\cO_{X^y}=
\begin{cases}\cO_{X_x^y} & \text{if } y\le x,\\
0 &\text{otherwise.}\\
\end{cases}
$$
Together with Proposition \ref{vanishing1}, it follows that
$$
\chi([\cO_{X_x}]\cdot[\cO_{X^y}])=
\begin{cases}1& \text{if }y\le x,\\
0& \text{otherwise}.\\
\end{cases}
$$
On the other hand, we have in $K(G/P)$:
$$
[\cO_{X^y}]=\sum_{z\in W^P,\;z\ge y} [\cO_{X^z}(-\partial X^z)]
$$
(more generally, for any union $Z$ of opposite Schubert varieties, we
have 
$[\cO_Z]=\sum_{z\in W^P,\;X^z\subseteq Z} [\cO_{X^z}(-\partial X^z)]$
by an easy induction, using the fact that intersections of unions of
opposite Schubert varieties are reduced.) It follows that
$$
\chi([\cO_{X_x}]\cdot[\cO_{X^y}(-\partial X^y)])=\delta_{x,y}.
$$
Thus, the $[\cO_{X_x}]$, $x\in W^P$ form a basis for $K(G/P)$;
further, the bilinear form 
$K(G/P)\times K(G/P)\longrightarrow\mZ,
~(u,v)\longmapsto \chi(u\cdot v)$ 
is non--degenerate, and the dual basis of the $[\cO_{X_x}]$ with
respect to this pairing consists of the $[\cO_{X^x}(-\partial X^x)]$.

It follows that a given class $u\in K(G/P\times G/P)$ is zero if and
only if 
$\chi(u\cdot [\cO_{X^y}(-\partial X^y)\otimes\cO_{X_z}])=0$
for all $y,z\in W^P$. Further,
$$
\chi([\cO_{\diag(G/P)}]\cdot[\cO_{X^y}(-\partial X^y)\otimes\cO_{X_z}])
=\chi([\cO_{X^y}(-\partial X^y)]\cdot[\cO_{X_z}])
=\delta_{y,z},
$$
whereas
$$\displaylines{
\chi(\sum_{x\in W^P} [\cO_{X_x}\otimes\cO_{X^x}(-\partial X^x)]
\cdot[\cO_{X^y}(-\partial X^y)\otimes\cO_{X_z}])=
\hfill\cr\hfill
=\sum_{x\in W^P}\chi([\cO_{X_x}]\cdot[\cO_{X^y}(-\partial X^y)])\;
\chi([\cO_{X^x}(-\partial X^x)]\cdot[\cO_{X_z}])
=\sum_{x\in W^P}\delta_{x,y}\delta_{x,z}=\delta_{y,z}.
\cr}$$
This completes the proof of $(*)$, and hence of the lemma.
\end{proof}

We now construct a degeneration of the diagonal of any Richardson
variety. Let $\theta:\mG_m\longrightarrow T$ be a regular dominant
one--parameter subgroup. Let $\cX$ be the closure in 
$G/P\times G/P\times \mA^1$ of the subset
$$
\{(x,\theta(s)x,s)~\vert~ x\in G/P,~s\in k^*\}.
$$
The variety $\cX$ is invariant under the action of $\mG_m\times T$
defined by
$$
(s,t) (x,y,z)= (t x,\theta(s)t y,sz).
$$
Consider the projections
$$
p_1,p_2:\cX\longrightarrow G/P, ~\pi:\cX\longrightarrow\mA^1.
$$
Clearly, $\pi$ is proper and flat, and its fibers identify with closed 
subschemes of $G/P\times G/P$ via $p_1\times p_2$; this identifies the
fiber at $1$ with $diag(G/P)\simeq G/P$. By equivariance,
every ``general'' fiber $\pi^{-1}(z)$, where $z\ne 0$, is also
isomorphic to $G/P$.

We shall denote the ``special'' (scheme--theoretical) fiber
$\pi^{-1}(0)$ by $F$, with projections
$$
q_1,q_2:F\longrightarrow G/P.
$$

Next let $v,w$ in $W^P$ such that $v\leq w$. Let $\cX_w^v$ be the
closure in $G/P\times G/P\times \mA^1$ of the subset
$$
\{(x,\theta(s)x,s)~\vert~ x\in X_w^v,~s\in k^*\}.
$$
This is a subvariety of $\cX\cap(X_w^v\times X_w^v\times\mA^1)$,
invariant under the action of $\mG_m\times T$. We shall denote the
restrictions of $p_1$, $p_2$, $\pi$ to $\cX_w^v$ by the same letters;
then $\pi$ is again proper and flat, and its ``general'' fibers are
isomorphic to $X_w^v$. Let $F_w^v$ be the ``special'' fiber, with
projections $q_1$, $q_2$ to $X_w^v$.

\begin{lemma}\label{flat}

\begin{enumerate}

\item The schemes $F$ and $F_w^v$ are reduced. Further,
$$
F=\bigcup_{x\in W^P} X_x \times X^x \text{ and }
F_w^v=\bigcup_{x\in W^P,\;v\leq x\leq w} X_x^v \times X_w^x.
$$

\item Choose a total ordering $\le_t$ of $W^P$ such that $x\le_t y$
whenever $x\le y$. For $x\in W^P$, let $\cO_{F,\le_t x}$ 
(resp.~$\cO_{F,\ge_t x}$) be the subsheaf of $\cO_F$ consisting of
those sections that vanish identically on $X_y\times X^y$ for each 
$y>_t x$ (resp.~$y<_t x$). Then the $\cO_{F,\le_t x}$
(resp.~$\cO_{F,\ge_t x})$ are an ascending (resp.~descending)
filtration of $\cO_F$, with associated graded
$$
\bigoplus_{x\in W^P} 
\cO_{X_x}\otimes\cO_{X^x}(-\partial X^x), \text{ resp. }
\bigoplus_{x\in W^P} 
\cO_{X_x}(-\partial X_x)\otimes\cO_{X^x}.
$$
The induced filtrations on the structure sheaf $\cO_{F_w^v}$ have
associated graded
$$
\bigoplus_{x\in W^P,\;v\leq x\leq w} 
\cO_{X_x^v}\otimes\cO_{X_w^x}(-(\partial X^x)_w), \text{ resp. }
\bigoplus_{x\in W^P,\;v\leq x\leq w} 
\cO_{X_x^v}(-(\partial X_x)^v)\otimes\cO_{X_w^x}.
$$
The induced map 
$$
\cO_{F_w^v,\le_t}\longrightarrow \gr_{\le t}\cO_{F_w^v}
$$ 
is just the restriction to $X_x^v\times X_w^x$; the same holds for the
induced map
$$
\cO_{F_w^v,\ge_t}\longrightarrow \gr_{\ge t}\cO_{F_w^v}.
$$
\end{enumerate}
\end{lemma}

\begin{proof}
(1) Let $x\in W^P$. We claim that
$$
C_x\times C^x\subseteq F.
$$
To check this, consider the subset $xC^1$ of $G/P$. This is an open
$T$--stable neighborhood of $e_x$ in $G/P$, isomorphic to affine space
where $T$ acts linearly with weights the $\alpha\in x(R^- - R_P)$. 
Choose corresponding coordinate functions $z_\alpha$ on $xC^1$,
then $C_x$ (resp.~$C^x$) is the closed subset of $xC^1$ where
$z_\alpha=0$ whenever $\alpha\in R^-$ (resp.~$\alpha\in R^+$).
Let $z=(z_\alpha)\in xC^1$, then 
$$
\theta(s)z=(s^{\langle\alpha,\theta\rangle}z_\alpha).
$$
Denote by $z_+$ (resp.~$z_-$) the point of $C_x$ (resp.~$C^x)$ with
coordinates $z_\alpha$, $\alpha\in R^+$ (resp.~$\alpha\in R^-$). 
Let $z'(s)$ be the point of $xC^1$ with $\alpha$--coordinate
$z_\alpha$ if $\alpha\in R^+$, and $\theta(s^{-1})z_\alpha$
otherwise. Since $\theta$ is regular dominant, we obtain
$$
\lim_{s\to 0} (z'(s),\theta(s)z'(s),s) = (z_+,z_-,0).
$$
And since $z_+$ (resp.~$z_-$) is an arbitrary point of $C_x$
(resp.~$C^x$), this proves our claim.

The claim implies that $F$ contains $\bigcup_{x\in W^P} X_x\times X^x$
as a reduced closed subscheme. Let $\cI$ be the ideal sheaf of this
closed subscheme in $\cO_F$; we regard $\cI$ as a coherent sheaf on
$G/P\times G/P$. Then we have in $K(G/P\times G/P)$:
$$
[\cI]=[\cO_F]-[\cO_{\bigcup_{x\in W^P} X_x\times X^x}]
=[\cO_{\diag(G/P}]-[\cO_{\bigcup_{x\in W^P} X_x\times X^x}] =0,
$$
where the first equality follows from the definition of $\cI$, the
second one from the fact that $\pi:\cX\to\mA^1$ is flat with fibers $F$
and $\diag(G/P)$, and the third one from Lemma \ref{formula}.
As a consequence, $\cI$ is trivial (e.g., since its Hilbert polynomial
is zero); this completes the proof for $F$.

In the case of $F_w^v$, notice that
$$
F_w^v\subseteq F\cap (X^v\times X_w) = 
(\bigcup_{x\in W^P} X_x\times X^x)\cap (X^v\times X_w)
=\bigcup_{x\in W^P, v\le x\le w} X_x^v\times X_w^x
$$
as schemes, since all involved scheme--theoretic intersections are
reduced. Further, we have in the Chow ring of $G/P\times G/P$:
$$\displaylines{
[F_w^v]=[\diag(X_w^v)]=[\diag(G/P)\cap(X^v\times X_w)]
=[\diag(G/P)]\cdot[X^v\times X_w]
\hfill\cr\hfill
=[F]\cdot [X^v\times X_w]
=[F\cap(X^v\times X_w)]
=\sum_{x\in W^P,\;v\le x\le w} [X_x^v\times X_w^x],
\cr}$$
since all involved intersections are proper and reduced.
It follows that $F_w^v$ equals 
$\bigcup_{x\in W^P,\;v\le x\le w} X_x^v\times X_w^x$.

(2) has been established in the case of $F$, at the beginning of the
proof of Lemma \ref{formula}. The general case is similar.
\end{proof}

Next let $\lambda$ be a dominant character of $P$. This yields
$T$--linearized line bundles $q_2^*L_\lambda$ on $F$ and on $F_w^v$,
together with ``adjunction'' maps 
$H^0(G/P,L_\lambda)\longrightarrow H^0(F,q_2^*L_\lambda)$ and  
$H^0(X_w^v,L_\lambda)\longrightarrow H^0(F_w^v,q_2^*L_\lambda)$.

\begin{proposition}\label{filtration1}
\begin{enumerate}

\item
These maps are isomorphisms, and the restriction map
$$
H^0(F,q_2^*L_\lambda)\longrightarrow H^0(F_w^v,q_2^*L_\lambda)
$$
is surjective. 

\item
The ascending filtration of $\cO_F$ yields an
ascending filtration of the $T$--module $H^0(F,q_2^*L_\lambda)$,
with associated graded 
$$
\bigoplus_{x\in W^P} H^0(X^x,L_\lambda(-\partial X^x)).
$$

\item
The image of this filtration under restriction to $F_w^v$ 
has associated graded 
$$
\bigoplus_{x\in W^P,\;v\leq x\leq w} 
H^0(X_w^x,L_\lambda(-(\partial X^x)_w)).
$$
Hence this is the associated graded of an ascending filtration 
$$
H^0(X_w^v,L_\lambda)_{\le_t x}, ~v\le x\le w
$$
of $H^0(X_w^v,L_\lambda)$, compatible with the $T$--action and with
restrictions to smaller Richardson varieties. 

\item
The subspace 
$$
H^0(X_w^v,L_\lambda)_{\le_t x}\subseteq
H^0(X_w^v,L_\lambda)
$$ 
consists of those sections that vanish identically on $X_w^y$ for all
$y>_t x$. Further, the map 
$$
H^0(X_w^v,L_\lambda)_{\le_t x} \longrightarrow 
\gr_x H^0(X_w^v,L_\lambda) = 
H^0(X_w^x,L_\lambda(-(\partial X^x)_w))
$$
is just the restriction to $X_w^x$.
\end{enumerate}
\end{proposition}

\begin{proof}
(1) We have
$$
H^0(F,q_2^*L_\lambda) = H^0(G/P,q_{2*}q_2^*L_\lambda)
= H^0(G/P,L_\lambda\otimes q_{2*}\cO_F)
$$
by the projection formula. Further, the associated graded of the
descending filtration of $\cO_F$ is acyclic
for $q_{2*}$; indeed, $H^i(X_x,\cO_{X_x}(-\partial X_x))=0$ for all
$i\geq 1$ and all $x\in W^P$, by Proposition \ref{vanishing1}. Notice
also that $H^0(X_x,\cO_{X_x}(-\partial X_x))=0$ for all $x\neq 1$,
since $\partial X_x$ is a nonempty subscheme of the complete variety
$X_x$. It follows that the natural map 
$\cO_{G/P}\longrightarrow q_{2*}\cO_F$ 
is an isomorphism. Hence the same holds for the map
$H^0(G/P,L_\lambda)\longrightarrow H^0(F,q_2^*L_\lambda)$. 

Likewise, the map 
$H^0(X_w^v,L_\lambda)\longrightarrow H^0(F_w^v,q_2^*L_\lambda)$
is an isomorphism as well. Since the restriction map
$H^0(G/P,L_\lambda)\longrightarrow H^0(X_w^v,L_\lambda)$ is surjective
by Proposition \ref{vanishing1}, the same holds for 
$H^0(F,q_2^*L_\lambda)\longrightarrow H^0(F_w^v,q_2^*L_\lambda)$.

(2) By Lemma \ref{flat} again, the ascending filtration of $\cO_F$
yields one on $q_2^*L_\lambda$, with associated graded 
$$
\bigoplus_{x\in W^P} \cO_{X_x}\otimes 
L_\lambda\vert_{X^x}(-\partial X^x).
$$
The latter is acyclic by Proposition \ref{vanishing1}. It follows that 
$H^0(F, q_2^*L_\lambda)$ has an ascending filtration with associated
graded as claimed. 

(3) is checked similarly. 

(4) We have
$$\displaylines{
H^0(X_w^v,L_\lambda)_{\le_t x} = H^0(F_w^v,q_2^*L_\lambda)_{\le_t x}
= H^0(F_w^v, q_2^*L_\lambda\otimes\cI_{\cup_{y>_t x} X_y^v\times X_w^y})
\hfill\cr\hfill
= H^0(X_w^v,L_\lambda\otimes q_{2*}\cI_{\cup_{y>_t x} X_y^v\times X_w^y})
\cr}$$
by the projection formula. Further, $q_{2*}\cO_{F_w^v}=\cO_{X_w^v}$
as seen in the proof of (1). It follows that 
$$
q_{2*}\cI_{\cup_{y>_t x} X_y^v\times X_w^y}=\cI_{\cup_{y>_t x} X_w^y}.
$$
This implies our statement.
\end{proof}

We now construct a similar filtration of the $T$--submodule
$$
H^0(X_w^v,L_\lambda(-(\partial_\lambda X_w)^v))\subseteq 
H^0(X_w^v,L_\lambda).
$$ 
For this, we define a sheaf on $F_w^v$ by
$$
q_2^*L_\lambda(-(\partial_\lambda X_w)^v)=
(q_2^*L_\lambda)\otimes_{\cO_{F_w^v}} 
\cI_{q_2^{-1}((\partial_\lambda X_w)^v)}.
$$
This is a subsheaf of $q_2^*L_\lambda$; it may differ from the
pullback sheaf of $L_\lambda(-(\partial_\lambda X_w)^v)$ under $q_2$.
We also have an ``adjunction'' map
$$
H^0(X_w^v,L_\lambda(-(\partial_\lambda X_w)^v))\longrightarrow 
H^0(F_w^v,q_2^*L_\lambda(-(\partial_\lambda X_w)^v)).
$$
In particular, we obtain a map
$$
H^0(X_w,L_\lambda(-\partial_\lambda X_w))\longrightarrow 
H^0(F_w,q_2^*L_\lambda(-\partial_\lambda X_w)).
$$

\begin{proposition}\label{filtration2}
\begin{enumerate}

\item
These maps are isomorphisms, and the restriction map 
$$
H^0(F_w,q_2^*L_\lambda(-\partial_\lambda X_w))\longrightarrow 
H^0(F_w^v,q_2^*L_\lambda(-(\partial_\lambda X_w)^v))
$$
is surjective. 

\item
The ascending filtration of $\cO_{F_w}$ yields
an ascending filtration of the $T$--module 
$H^0(F_w,q_2^*L_\lambda(-\partial_\lambda X_w))$,
with associated graded
$$
\bigoplus_{x\in W^P,\;x \leq w} 
H^0(X_w^x,L_\lambda(-(\partial_\lambda X_w)^x - (\partial X^x)_w)).
$$

\item
The image of this filtration under restriction to $F_w^v$
has associated graded
$$
\bigoplus_{x\in W^P,\;v\leq x\leq w} 
H^0(X_w^x,L_\lambda(-(\partial_\lambda X_w)^x - (\partial X^x)_w)).
$$
Hence this is also the associated graded of an ascending filtration 
$$
H^0(X_w^v,L_\lambda(-(\partial_{\lambda} X_w)^v))_{\le_t x},~v\le x\le w
$$
of $H^0(X_w^v,L_\lambda(-(\partial_{\lambda} X_w)^v))$, compatible
with the $T$--action and with restrictions to smaller Richardson
varieties. 

\item
The subspace 
$$
H^0(X_w^v,L_\lambda(-(\partial_\lambda X_w)^v))_{\le_t x}
\subseteq H^0(X_w^v,L_\lambda)
$$ 
consists of those sections that vanish identically on
$(\partial_\lambda X_w)^v$ and on $X_w^y$ for all
$y>_t x$. Further, the map 
$$\displaylines{
H^0(X_w^v,L_\lambda(-(\partial_\lambda X_w)^v))_{\le_t x} 
\longrightarrow 
\gr_x H^0(X_w^v,L_\lambda(-(\partial_\lambda X_w)^v)) = 
\hfill\cr\hfill
= H^0(X_w^x,L_\lambda(-(\partial_\lambda X_w)^x-(\partial X^x)_w))
\cr}$$
is just the restriction to $X_w^x$.
\end{enumerate}
\end{proposition}

\begin{proof}
It follows from Lemma \ref{flat} that the sheaf
$q_2^*L_\lambda(-(\partial_\lambda X_w)^v)$ on $F_w^v$ admits a
descending filtration with associated graded
$$
\bigoplus_{x\in W^P,v\leq x\leq w} \cO_{X_x^v}(-(\partial X_x)^v)
\otimes L_\lambda\vert_{X_w^x}(-(\partial_{\lambda}X_w)^v),
$$
and an ascending filtration with associated graded
$$
\bigoplus_{x\in W^P,v\leq x\leq w} \cO_{X_x^v}\otimes 
L_\lambda\vert_{X_w^x}(-(\partial_{\lambda}X_w)^v-(\partial X_x)^v).
$$
As in the proof of Proposition \ref{filtration1}, the associated
graded of the first filtration is acyclic for $q_{2*}$; it follows
that the adjunction map is an isomorphism. Further, the restriction map
$$
H^0(X_w,L_\lambda(-\partial_\lambda X_w))\longrightarrow 
H^0(X_w^v,L_\lambda(-(\partial_\lambda X_w)^v))
$$
is surjective, by Corollary \ref{surj}. Finally, the associated graded
of the second filtration is acyclic, by Proposition
\ref{vanishing2}. These facts imply our statements, as in the proof of
Proposition \ref{filtration1}.
\end{proof}

\medskip

\noindent
{\bf Remarks.}
\begin{enumerate}

\item
By Proposition \ref{filtration1}, the $H^0(G/P,L_\lambda)_{\le_t x}$
are $B^-$--submodules of $H^0(G/P,L_\lambda)$. Likewise, the descending
filtration of $\cO_F$ yields a descending filtration of
$H^0(G/P,L_\lambda)$ by $B$--submodules $H^0(G/P,L_\lambda)_{\ge_t x}$,
consisting of those sections that vanish on $X_y$ whenever $y<_t x$.

\item
We may have defined directly the preceding filtrations by
Propositions \ref{filtration1} (4)  and \ref{filtration2} (4), 
without using the degeneration of the diagonal constructed in Lemma 
\ref{flat}. In fact, this alternative definition suffices for the
construction of a standard basis in the next section. But the
degeneration of the diagonal will play an essential r\^ole in the
section on standard products.

\end{enumerate}
 
\section{Construction of a standard basis}

In this section, we fix a dominant weight $\lambda$ and we consider
Richardson varieties in $G/P$, where $P=P_\lambda$. We shall construct
a basis of $H^0(G/P,L_\lambda)$ adapted to the filtrations of
Propositions \ref{filtration1} and \ref{filtration2}. We first prove
the key

\begin{lemma}\label{lift}
Let $v\le w\in W^\lambda$. Then any element of 
$H^0(X_w^v,L_\lambda(-(\partial X_w)^v-(\partial X^v)_w))$
can be lifted to an element of $H^0(G/P,L_\lambda)$ that vanishes
identically on all Schubert varieties $X_y$, $y\not\ge w$, and on all
opposite Schubert varieties $X^x$, $x\not\le v$. 
\end{lemma}

\begin{proof}
Put 
$$
X=X_w^v\text{ and }
Y=(\bigcup_{y\not\ge w} X_y)\cup(\bigcup_{x\not\le v} X^x).
$$
Notice that 
$$
X\cap Y = (\partial X_w)^v \cup (\partial X^v)_w
$$
(as schemes), since any intersection of unions of Richardson varieties
is reduced. This yields an exact sequence
$$
0 \longrightarrow \cI_{X\cup Y} \longrightarrow \cI_Y
\longrightarrow \cI_Y\otimes_{\cO_{G/P}}\cO_X\simeq\cO_X(-X\cap Y) 
\longrightarrow 0.
$$
Tensoring by $L_\lambda$ and taking the associated long exact sequence
of cohomology groups yields an exact sequence
$$
H^0(G/P,L_\lambda\otimes \cI_Y) \longrightarrow 
H^0(X,L_\lambda(-X\cap Y)) \longrightarrow
H^1(X,L_\lambda\otimes \cI_{X\cup Y}).
$$
Further, $H^1(X,L_\lambda\otimes \cI_{X\cup Y})=0$ by Lemma
\ref{vanishing}; this completes the proof.
\end{proof}

\begin{defn}\label{construction}
For any $v\le w\in W^\lambda$, let
$$
H_w^v(\lambda)=
H^0(X_w^v,L_\lambda(-(\partial X_w)^v -(\partial X^v)_w))
$$
and
$$
\chi_w^v(\lambda)=\{ \text{the weights of the }T--\text{module }
H_w^v(\lambda)\},
$$
these weights being counted with multiplicity. Let
$$
\{p_{w,v}^\xi,~\xi\in \chi_w^v(\lambda)\}
$$
be a basis for $H_w^v(\lambda)$, where each $p_{w,v}^\xi$ is a
$T$--eigenvector of weight $\xi$. 

For any triple $(w,v,\xi)$ as above, let $p_\pi$ be a lift of
$p_{w,v}^\xi$ in $H^0(G/P,L_\lambda)$ such that:

\noindent
$p_\pi$ is a $T$--eigenvector of weight $\xi$, and 

\noindent
$p_\pi$ vanishes identically on all $X_y$, $y\not\ge w$ and on all
$X^x$, $x\not\le v$.
\end{defn}

\noindent
(The existence of such lifts follows from Lemma \ref{lift}.) If $v=w$,
then $X_w^v$ consists of the point $e_w$, and hence
$\chi_w^v(\lambda)$ consists of the weight $-w(\lambda)$. We then
denote the unique $p_{w,v}^\xi$ by $p_w$. Its lift to
$H^0(G/P,L_\lambda)$ is unique; it is the extremal weight vector $p_w$.

\begin{defn}
Let $\pi=(w,v,\xi)$ be as in Definition \ref{construction}. We set
$i(\pi)=w$, $e(\pi)=v$, and call them respectively the initial and end
elements of $\pi$.
\end{defn}

By construction of the $p_\pi$ and Lemma \ref{richardson}, we obtain:

\begin{lemma}\label{nonvanishing1}
With notations as above, we have for $x, y\in W^\lambda$:
$$
p_\pi\vert_{X_y^x}\not= 0\Longleftrightarrow 
X_{i(\pi)}^{e(\pi)}\subseteq X_y^x
\Longleftrightarrow x \le e(\pi)\le i(\pi)\le y.
$$
\end{lemma}

\begin{proposition}\label{bound}
The restrictions to $X^w_v$ of the $p_\pi$ where $i(\pi)=w$, 
$e(\pi)\ge v$ form a basis for the $T$--module
$H^0(X_w^v,L_\lambda(-(\partial X_w)^v))$, adapted to its ascending
filtration $\le_t$ of Proposition \ref{filtration2}.
\end{proposition}

\begin{proof}
By construction, $p_\pi$ vanishes identically on $X^x$ for any
$x\not\le e(\pi)$, and hence for any $x>_t e(\pi)$. Thus, 
$p_\pi\in H^0(G/P,L_\lambda)_{\le_t e(\pi)}$ by Proposition
\ref{filtration1}. Further, the image of $p_\pi$ in the associated
graded is just its restriction to $X^{e(\pi)}$. 

Together with Lemma \ref{nonvanishing1}, it follows that the
restrictions  of the $p_\pi$ to $X_w^v$ belong to  
$H^0(X_w^v,L_\lambda(-(\partial X_w)^v))_{\le_t e(\pi)}$, 
and that their images in the associated graded $H_w^{e(\pi)}(\lambda)$
are the restrictions of the $p_\pi$ to $X_w^{e(\pi)}$; by
construction, these images form a basis of $H_w^{e(\pi)}(\lambda)$.
\end{proof}

Now the $T$--module $H^0(X_w^v,L_\lambda)$
has a descending filtration by the submodules
$$
H^0(X_w^v,L_\lambda(-(\partial X_w)^v))_{\ge_t x}
$$
consisting of those sections that vanish identically on $X_y^v$
whenever $y<_t x$. And like in Proposition \ref{filtration1}, the
associated graded is
$$
\bigoplus_{x\in W^P, v\le x\le w} 
H^0(X_x^v,L_\lambda(-(\partial X_x)^v)).
$$
Further, we may check as in the proof of Proposition \ref{bound} that
$$
p_\pi\vert_{X_w^v} \in 
H^0(X_w^v,L_\lambda(-(\partial X_w)^v))_{\ge_t i(\pi)}
$$
whenever $i(\pi)\ge w$, and the image of $p_\pi$ in the
associated graded is just its restriction to $X_{i(\pi)}^v$. Together
with Proposition \ref{bound}, this implies 

\begin{proposition}\label{adapted}
The restrictions to $X_w^v$ of the $p_\pi$ where 
$v\le e(\pi)\le i(\pi)\le w$ form a basis of $H^0(X_w^v,L_\lambda)$;
the $p_\pi$ where $v\not\le e(\pi)$ or $i(\pi)\not\le w$ form a basis
of the kernel of the restriction map
$H^0(G/P,L_\lambda)\longrightarrow H^0(X_w^v,L_\lambda)$. 
\end{proposition}

In view of Proposition \ref{adapted}, the restriction to $p_\pi$ 
to $X_y^x$, where $x\le e(\pi)\le i(\pi) \le y$, will be denoted 
by just $p_\pi$.

\begin{defn} 
Set
$$
\Pi(\lambda)=\{(v,w,\xi)~\vert~
v,w\in W^\lambda,\;v\le w,\;\xi\in\chi_w^v(\lambda)\}. 
$$
For any $v, w\in W^\lambda, v\le w$, set
$$
\Pi_w^v(\lambda):=\{\pi~\vert~v \le e(\pi)\le i(\pi)\le w\}.
$$ 
In view of Lemma \ref{nonvanishing1}, we have,
$\Pi_w^v(\lambda)=\{\pi\in\Pi(\lambda)~\vert~ p_\pi\vert_{X_w^v}\neq 0\}$.

More generally, for a union $Z$ of Richardson varieties, define
$$
\Pi_Z(\lambda)=\{\pi\in\Pi(\lambda)~\vert~ p_\pi\vert_Z\neq 0\}.
$$
\end{defn}

\begin{theorem}\label{basis1}
Let $Z$ be a union of Richardson varieties. Then
$\{p_\pi\vert_Z,\pi \in \Pi_Z(\lambda)\}$ is a basis for
$H^0(Z,L_\lambda)$, and 
$\{p_\pi,\pi \in \Pi(\lambda) -\Pi_Z(\lambda)\}$ 
is a basis for the kernel of the restriction map 
$H^0(G/P,L_\lambda)\longrightarrow H^0(Z,L_\lambda)$.
\end{theorem}

\begin{proof}
By definition of $\Pi_Z(\lambda)$, it suffices to prove the first
assertion. Let $Z=\cup_{i=1}^r\ X_{w_i}^{v_i}$. We shall prove the
result by induction on $r$ and dim $Z$. Write $Z=X\cup Y$ where
$X=X_{w_i}^{v_i}$ for some $i$, and $\dim X= \dim Z$. Then $X\cap Y$
is a union of Richardson varieties of dimension $<\dim Z$.
Consider the exact sequence
\begin{equation*}
0\rightarrow \cO_Z=\cO_{X\cup Y}\rightarrow \cO_X \oplus
\cO_Y\rightarrow \cO_{X\cap Y}\rightarrow 0.\tag{*}
\end{equation*}
Tensoring by $L_\lambda$, taking global sections and using the
vanishing of $H^1(Z,L_\lambda)$ (Lemma \ref{vanishing}), we obtain the
exact sequence
$$
0 \longrightarrow H^0(Z, L_\lambda) \longrightarrow 
H^0(X,L_\lambda)\oplus H^0(Y, L_\lambda) \longrightarrow 
H^0(X\cap Y, L_\lambda) \longrightarrow 0.
$$ 
In particular, denoting $\dim H^0(Z, L_\lambda)$ by 
$h^0(Z, L_\lambda)$ etc., we obtain, 
$$
h^0(Z, L_\lambda)=h^0(X, L_\lambda)+h^0(Y, L_\lambda)
-h^0(X\cap Y, L_\lambda).
$$ 
We have by hypothesis (and induction hypothesis),
$h^0(X, L_\lambda)=\#\Pi_X(\lambda)$, 
$h^0(Y,L_\lambda)=\#\Pi_Y(\lambda)$, 
$h^0(X\cap Y, L_\lambda)=\#\Pi_{X\cap Y }(\lambda)$. 
Thus we obtain,
\begin{equation*}
h^0(Z, L_\lambda)=\#\Pi_X(\lambda)+\#\Pi_Y(\lambda)-
\#\Pi_{X\cap Y}(\lambda).\tag{1}
\end{equation*}
On the other hand we have,
\begin{equation*}
\Pi_Z(\lambda)=(\Pi_X(\lambda)\,\dot\cup\,\Pi_Y(\lambda))
\setminus\Pi_{X\cap Y}(\lambda).\tag{2}
\end{equation*}
From (1) and (2), we obtain, $h^0(Z, L_\lambda)=\#\Pi_Z(\lambda)$.
Further, the $p_\pi\vert_Z$, $\pi\in\Pi_Z(\lambda)$, span
$H^0(Z,L_\lambda)$ (since the $p_\pi$, $\pi\in\Pi(\lambda)$, span
$H^0(G/P,L_\lambda)$, and the restriction map 
$H^0(G/P,L_\lambda)\to H^0(Z,L_\lambda)$ is surjective). Thus, the 
$p_\pi\vert_Z$, $\pi\in\Pi_Z(\lambda)$, are a basis of
$H^0(Z,L_\lambda)$.
\end{proof}

\section{Standard monomials}

Let $\lambda$, $\mu$ be dominant weights such that 
$P_\lambda = P_\mu:=P$. Consider the product map
$$
H^0(G/P,L_\lambda)\otimes H^0(G/P,L_\mu)\longrightarrow 
H^0(G/P,L_{\lambda+\mu}).
$$
This map is surjective by \cite{Ram2} 2.2 and 3.5. Using Proposition
\ref{vanishing1}, it follows that the product map
$$
H^0(X_w^v,L_\lambda)\otimes H^0(X_w^v,L_\mu)\longrightarrow 
H^0(X_w^v,L_{\lambda+\mu})
$$
is also surjective, for any $v\leq w$ in $W^P$. We shall construct 
a basis for $H^0(X_w^v,L_{\lambda+\mu})$ from the bases of
$H^0(X_w^v,L_\lambda)$, $H^0(X_w^v,L_\mu)$ obtained in Theorem
\ref{basis1}. For this, we need the following 

\begin{defn}
Let $v,w \in W^P, v\le w$. Let $\varphi\in\Pi_w^v(\lambda)$ and
$\psi\in\Pi_w^v(\mu)$. The pair $(\varphi,\psi)$ is called
standard on $X_w^v$ if
$$
v\le e(\psi) \le i(\psi) \le e(\varphi) \le i(\varphi) \leq w.
$$ 
Then the product $p_{\varphi}p_{\psi}\in H^0(G/P,L_{\lambda+\mu})$ 
is called standard on $X_w^v$ as well.
\end{defn}

Clearly, we have

\begin{lemma}\label{nonvanishing2}
Let $p_\varphi p_\psi$ be a standard product on $G/P$ and let 
$v\le w\in W^P$. Then 
$$
p_\varphi p_\psi \vert_{X_w^v}\neq 0 \Longleftrightarrow
v\le e(\psi)\le i(\varphi)\le w.
$$
\end{lemma}

\begin{proposition}\label{products}
The standard products on $X_w^v$ form a basis of
$H^0(X_w^v,L_{\lambda+\mu})$. The standard products on $G/P$ that are
not standard on $X_w^v$ form a basis of the kernel of the restriction map 
$H^0(G/P,L_{\lambda+\mu})\longrightarrow H^0(X_w^v,L_{\lambda+\mu})$.
\end{proposition}

\begin{proof}
Consider the $T$--linearized invertible sheaf
$q_1^*L_\mu\otimes q_2^*L_\lambda$ on $F_w^v$. By Lemma \ref{flat},
the ascending filtration of $\cO_{F_w^v}$ yields one of that sheaf,
with associated graded
$$
\bigoplus_{x\in W^P,v\leq x\leq w} 
L_\mu\vert_{X_x^v}(-(\partial X_x)^v)
\otimes L_\lambda\vert_{X_w^x}.
$$
By Proposition \ref{vanishing1}, the latter sheaf is acyclic. This
yields an ascending filtration of the $T$--module  
$H^0(F_w^v,q_1^*L_\mu\otimes q_2^*L_\lambda)$, with associated graded
$$
\bigoplus_{x\in W^P,v\leq x\leq w}
H^0(X_x^v,L_\mu(-(\partial X_x)^v))\otimes
H^0(X_w^x,L_\lambda);
$$
it also follows that $H^i(F_w^v,q_1^*L_\mu\otimes q_2^*L_\lambda)=0$ 
for all $i\geq 1$.

By Proposition \ref{filtration1}, we may identify
$H^0(X_w^v,L_\lambda)$ with $H^0(F_w^v,q_2^*L_\lambda)$; likewise, we
may identify $H^0(X_w^v,L_\mu)$ with $H^0(F_w^v,q_1^*L_\mu)$. Using
the multiplication map 
$$
H^0(F_w^v,q_1^*L_\mu) \otimes H^0(F_w^v,q_2^*L_\lambda)
\longrightarrow H^0(F_w^v, q_1^*L_\mu \otimes q_2^*L_\lambda),
$$
this defines ``dot products'' in 
$H^0(F_w^v, q_1^*L_\mu \otimes q_2^*L_\lambda)$.

Let $x\in W^P$ such that $v\leq x\leq w$. Recall that the $p_\psi$, 
$v\leq e(\psi)\leq i(\psi)=x$, are a basis of
$H^0(X_x^v,L_\mu(-(\partial X_x)^v))$. Further, the $p_\varphi$, 
$x\leq e(\varphi)\leq i(\varphi)\leq w$, are a basis
of $H^0(X_w^x,L_\lambda)$. Thus, the dot products 
$p_\psi \cdot p_\varphi$, where there exists $x\in W^P$ such that
$$
v\leq e(\psi)\leq i(\psi)=x \text{ and }
x\leq e(\varphi)\leq i(\varphi)\leq w,
$$ 
restrict to a basis of 
$H^0(X_x^v,L_\mu(-(\partial X_x)^v)) \otimes H^0(X_w^x,L_\lambda)$.
By construction of the filtration of 
$H^0(F_w^v, q_1^*L_\mu \otimes q_2^*L_\lambda)$, it follows that the
standard dot products are a basis of that space.

Consider now the $T$--linearized invertible sheaf 
$p_1^*L_\mu \otimes p_2^*L_\lambda$ on $\cX_w^v$. This sheaf is flat
on $\mA^1$; by vanishing of 
$H^1(F_w^v,q_1^*L_\mu \otimes q_2^*L_\lambda)$ and semicontinuity, it 
follows that the restriction 
$$
H^0(\cX_w^v,p_1^*L_\mu \otimes p_2^*L_\lambda)\longrightarrow
H^0(F_w^v,q_1^*L_\mu \otimes q_2^*L_\lambda)
$$
is surjective, and that 
$H^0(\cX_w^v,p_1^*L_\mu \otimes p_2^*L_\lambda)$ is a free module over
$H^0(\mA^1,\cO_{\mA_1})=k[z]$, generated by any lift of
its quotient space $H^0(F_w^v,q_1^*L_\mu \otimes q_2^*L_\lambda)$.

We now construct such a lift, as follows. Consider the adjunction maps
$$
H^0(X_w^v,L_\lambda)\longrightarrow H^0(\cX_w^v,p_2^*L_\lambda) 
\text{ and } 
H^0(X_w^v,L_\mu)\longrightarrow H^0(\cX_w^v,p_1^*L_\mu).
$$
These yield dot products $p_\psi \cdot p_\varphi$ in 
$H^0(\cX_w^v,p_1^*L_\mu\otimes p_2^*L_\lambda)$
which lift the corresponding products in 
$H^0(F_w^v,q_1^*L_\mu\otimes q_2^*L_\lambda)$.
Since the latter standard products are a basis of that space, the
standard dot products $p_\psi \cdot p_\varphi$ are a basis of 
$H^0(\cX_w^v,p_1^*L_\mu \otimes p_2^*L_\lambda)$ over
$k[z]$. Therefore, they restrict to a basis of the space 
of sections of $p_1^*L_\mu\otimes p_2^*L_\lambda$ over any fiber of
$\pi$. But the fiber at $1$ is $diag(X_w^v)$, and the
restriction of $p_1^*L_\mu\otimes p_2^*L_\lambda$ to that fiber is
just $L_{\lambda+\mu}$ whereas the restrictions of the dot products
are just the usual products. We have proved that the standard products
on $X_w^v$ form a basis of $H^0(X_w^v,L_{\lambda+\mu})$. 

To complete the proof, notice that any standard product on $G/P$ that
is not standard on $X_w^v$ vanishes identically on that subvariety, by
Lemma \ref{nonvanishing2}.
\end{proof}

\noindent
{\bf Remark.} The proof of Proposition \ref{products} relies on the
fact that the special fiber $F_w^v$ of the flat family
$\pi:\cX_w^v\to\mA^1$ equals 
$\bigcup_{x\in W^P,v\le x\le w} X_x^v\times X_w^x$. Conversely, 
this fact can be recovered from Proposition \ref{products}, as
follows. 

We have the equalities of Euler characteristics:
$$
\chi(F,q_1^*L_\mu\otimes q_2^*L_\lambda)=\chi(G/P,L_{\lambda+\mu})
=\sum_{x\in W^P} \chi(X_x,L_\mu(-\partial X_x))\;\chi(X^x,L_\lambda),
$$
where the first equality holds by flatness of $\pi$, and the second
one by Propositions \ref{bound}, \ref{adapted} and \ref{products}.
It follows that 
$$
\chi(F,q_1^*L_\mu\otimes q_2^*L_\lambda)=
\chi(\bigcup_{x\in W^P} X_x\times X^x,q_1^*L_\mu\otimes q_2^*L_\lambda).
$$ 
Since $F$ contains $\bigcup_{x\in W^P} X_x\times X^x$
by the first claim in the proof of Lemma \ref{flat}, and $\lambda$,
$\mu$ are arbitrary dominant $P$--regular weights, it follows that 
$F=\bigcup_{x\in W^P} X_x\times X^x$ (e.g., since both have the same
Hilbert polynomial). Now the argument of Lemma \ref{flat} yields 
$F_w^v=\bigcup_{x\in W^P,v\le x\le w} X_x^v\times X_w^x$.

\medskip

We now extend Proposition \ref{products} to unions of Richardson
varieties. 

\begin{defn}
Let $\Pi(\lambda,\mu)$ be the set of all standard pairs
$(\varphi,\psi)$ where $\varphi\in\Pi(\lambda)$ and 
$\psi\in\Pi(\mu)$. For $v\leq w\in W^P$, let $\Pi_w^v(\lambda,\mu)$ 
be the subset of standard pairs on $X_w^v$. In view of Lemma 
\ref{nonvanishing2}, we have
$$
\Pi_w^v(\lambda,\mu)=\{(\varphi,\psi)\in\Pi(\lambda,\mu)~\vert~
p_\varphi p_\psi\vert_{X_w^v}\neq 0\}.
$$
Finally, for a union $Z$ of Richardson varieties, let
$$
\Pi_Z(\lambda,\mu)=\{(\varphi,\psi)\in\Pi(\lambda,\mu)~\vert~
p_\varphi p_\psi\vert_Z\neq 0\}.
$$
\end{defn}

Now arguing as in the proof of Theorem \ref{basis1}, we obtain

\begin{theorem}\label{basis2}
Let $Z$ be a union of Richardson varieties in $G/P$. Then the products
$p_\varphi p_\psi$, where $(\varphi,\psi)\in\Pi_Z(\lambda,\mu)$, form a
basis of $H^0(Z,L_{\lambda+\mu})$. The products $p_\varphi p_\psi$,
where $(\varphi,\psi)\in\Pi(\lambda,\mu) - \Pi_Z(\lambda,\mu)$, form a
basis of the kernel of the restriction map
$H^0(G/P,L_{\lambda+\mu})\longrightarrow H^0(Z,L_{\lambda+\mu})$.
\end{theorem}

\begin{corollary}\label{relations}
For any $\varphi\in\Pi(\lambda)$ and $\psi\in\Pi(\mu)$, the product 
$p_\varphi p_\psi\in H^0(G/P,L_{\lambda+\mu})$ is a linear combination
of standard products $p_{\varphi'} p_{\psi'}$ where 
$i(\varphi')\ge i(\varphi)$ and $e(\psi')\le e(\psi)$.
\end{corollary}

\begin{proof}
Notice that $p_\varphi p_\psi$ vanishes identically on all $X_y$
where $y\not\ge i(\varphi)$, and on all $X^x$ where 
$x\not\le e(\psi)$. By Theorem \ref{basis2}, it follows that 
$p_\varphi p_\psi$ is a linear combination of standard products 
$p_{\varphi'} p_{\psi'}$, where $i(\varphi)'\not\le y$ whenever
$y\not\ge i(\varphi)$, and $e(\psi)'\not\ge x$ whenever
$x\not\le e(\psi)$. But this means exactly that 
$i(\varphi')\ge i(\varphi)$ and $e(\psi')\le e(\psi)$.
\end{proof}

Next we consider a family of dominant weights
$\lambda_1,\ldots,\lambda_m$ such that 
$P=P_{\lambda_1}=\cdots=P_{\lambda_m}$.  
For any union $Z$ of Richardson varieties in $G/P$, we shall construct
a basis of $H^0(Z, L_{\lambda_1+\cdots+\lambda_m})$, in terms of
{\it standard monomials of degree} $m$. These are defined as follows.

\begin{defn}
Let $\pi_i\in \Pi(\lambda_i)$ for $1\le i\le m$. Then the sequence 
${\underline{\pi}}:=(\pi_1,\pi_2,\ldots , \pi_m)$ is standard if 
$$
e(\pi_m)\le i(\pi_m)\le\cdots\le e(\pi_1)\le i(\pi_1).
$$
Further, let $v,w\in W^P$ such that $v\leq w$; then 
${\underline{\pi}}$ is standard on $X_w^v$ if 
$$
v\le e(\pi_m)\le i(\pi_m)\le\cdots\le e(\pi_1)\le i(\pi_1)\le w.
$$

Finally, ${\underline{\pi}}$ is standard on $Z=\cup X_{w_i}^{v_i}$ if
it is standard on $X_{w_i}^{v_i}$ for some $i$.
\end{defn}

Set
$$
\Pi_w^v(\lambda_1,\ldots,\lambda_m)=
\{\underline{\pi}=(\pi_1,\pi_2,\ldots, \pi_m)\,\vert \,
\underline{\pi} \text{ is standard on }X_w^v\},
$$
$$
\Pi_Z(\lambda_1,\ldots,\lambda_m)=
\{\underline{\pi}=(\pi_1,\pi_2,\ldots ,\pi_m)\,\vert \,
\underline{\pi} \text{ is standard on }Z\}.
$$
\begin{defn}
Given ${\underline{\pi}}=(\pi_1,\pi_2,\ldots , \pi_m)$, set
$p_{{\underline{\pi}}}:=p_{\pi_1}\cdots p_{\pi_m}$.

\noindent 
Note that 
$p_{{\underline{\pi}}}\in H^0(G/P,L_{\lambda_1+\cdots+\lambda_m})$. 
If ${\underline{\pi}}$ is standard, then we call $p_{{\underline{\pi}}}$
a standard monomial on $G/P$. If ${\underline{\pi}}$ is
standard on $X_w^v$ (resp. $Z$) , then we call
$p_{{\underline{\pi}}}$ a standard monomial on $X_w^v$ (resp. $Z$).
\end{defn}

By Theorem \ref{basis2} and induction on $m$, we obtain

\begin{corollary}\label{monomials}
Let $Z$ be a union of Richardson varieties in $G/P$ and let
$\lambda_1,\ldots,\lambda_m$ be dominant weights such that
$P=P_{\lambda_1}=\cdots =P_{\lambda_m}$. Then the monomials 
$p_{\underline{\pi}}$ where 
$\underline{\pi} \in \Pi_Z(\lambda_1,\ldots,\lambda_m)$ 
form a basis of $H^0(Z, L_{\lambda_1+\cdots+\lambda_m})$.
Further, the monomials $p_{\underline{\pi}}$ where 
$\underline{\pi} \in 
\Pi(\lambda_1,\ldots,\lambda_m)-\Pi_Z(\lambda_1,\ldots,\lambda_m)$, 
form a basis of the kernel of the restriction map  
$H^0(G/P,L_{\lambda_1+\cdots+\lambda_m})\longrightarrow 
H^0(Z, L_{\lambda_1+\cdots+\lambda_m})$.
\end{corollary}

As an application, we determine the equations of unions of Richardson
varieties in their projective embeddings given by very ample line
bundles on $G/P$. Let $\lambda$ be a dominant $P$--regular weight. For
any $\pi_1,\pi_2\in\Pi(\lambda)$, we have in $H^0(G/P,L_{2\lambda})$:
$$
p_{\pi_1}p_{\pi_2}-\sum a_{\pi'_1,\pi'_2} p_{\pi'_1} p_{\pi'_2}=0,
$$
where $a_{\pi'_1,\pi'_2}\in k$ and the sum is over those standard pairs
$(\pi'_1,\pi'_2)\in\Pi(\lambda,\lambda)$ such that 
$i(\pi'_1)\ge i(\pi_1)$ and $e(\pi'_2)\le e(\pi_2)$ (as follows from
Corollary \ref{relations}). 

\begin{defn}
The preceding elements 
$p_{\pi_1}p_{\pi_2}-\sum a_{\pi'_1,\pi'_2} \,p_{\pi'_1}\, p_{\pi'_2}$
when regarded in $S^2H^0(G/P,L_\lambda)$, will be called the quadratic
straightening relations.
\end{defn} 

\begin{corollary}\label{generators}
Let $\lambda$ be a regular dominant character of $P$.
\begin{enumerate}

\item
The multiplication map 
$$
\bigoplus_{m=0}^{\infty} S^m H^0(G/P,L_\lambda) \longrightarrow 
\bigoplus_{m=0}^{\infty} H^0(G/P,L_{m\lambda})
$$
is surjective, and its kernel is generated as an ideal by the
quadratic straightening relations.

\item
For any union $Z$ of Richardson varieties in $G/P$, the restriction map
$$
\bigoplus_{m=0}^{\infty} H^0(G/P,L_{m\lambda}) \longrightarrow 
\bigoplus_{m=0}^{\infty} H^0(Z,L_{m\lambda})
$$
is surjective. Its kernel is generated as an ideal by the $p_\pi$,
$\pi\in\Pi(\lambda)-\Pi_Z(\lambda)$ together with the standard products
$p_{\pi_1}p_{\pi_2}$ where $i(\pi_1)\not\le w$ or $e(\pi_2)\not\ge v$
whenever $X_w^v$ is an irreducible component of $Z$. If in addition
$Z$ is a union of Richardson varieties $X_w^v$ all having the same
$w$, then the $p_\pi$, $\pi\in\Pi(\lambda)-\Pi_Z(\lambda)$ suffice.
\end{enumerate}
\end{corollary}

\begin{proof}
(1) By \cite{Ram2} Theorem 3.11, the multiplication map is surjective,
and its kernel is generated as an ideal by the kernel $K$ of the map
$S^2H^0(G/P,L_\lambda)\longrightarrow H^0(G/P,L_{2\lambda})$. 
Let $J$ be the subspace of $S^2H^0(G/P,L_\lambda)$ generated by all
quadratic straightening relations. Then $J\subseteq K$, and the quotient
space $S^2H^0(G/P,L_\lambda)/J$ is spanned by the images of the
standard products. Further, their images in 
$S^2H^0(G/P,L_\lambda)/K\simeq H^0(G/P,L_{2\lambda})$ form a basis, by
Proposition \ref{products}. It follows that $J=K$.

(2) The first assertion follows from Lemma \ref{vanishing}. Consider a
standard monomial 
$p_{\underline{\pi}}=p_{\pi_1}\cdots p_{\pi_m}\in H^0(G/P,L_{m\lambda})$.
By Corollary \ref{monomials}, $p_{\underline{\pi}}$ vanishes
identically on $Z$ if and only if: $i(\pi_1)\not\le w$ or 
$e(\pi_m)\not\ge v$ for all irreducible components $X_w^v$. This
amounts to: $p_{\pi_1} p_{\pi_m}$ vanishes identically on $Z$. If in
addition $w$ is independent of the component, then $p_{\pi_1}$ or
$p_{\pi_m}$ vanishes identically on $Z$; further, $p_{\pi_1}p_{\pi_m}$
is a standard product on $G/P$. This implies the remaining
assertions, since the kernel of 
$H^0(G/P,L_{m\lambda})\longrightarrow H^0(Z,L_{m\lambda})$
is spanned by those standard monomials on $G/P$ that are not standard
on $Z$ (Corollary \ref{monomials}).
\end{proof}

\noindent
{\bf Remark.} In particular, the $p_\pi$, where
$\pi\in\Pi(\lambda)-\Pi_Z(\lambda)$, generate the homogeneous ideal of
$Z$ in $G/P$, whenever $Z$ is a union of Schubert varieties (or a union
of opposite Schubert varieties). But this does not extend to arbitrary
unions of Richardson varieties, as shown by the obvious example where
$G/P=\mP^1$, $Z=\{0,\infty\}$ and $L_\lambda=\cO(1)$; then
$\Pi(\lambda)=\Pi_Z(\lambda)$.

\section{Weights of classical type}

In this section, we shall determine the ``building blocks''
$$
H_w^v(\lambda)=H^0(X_w^v,
L_\lambda(-(\partial X_w)^v-(\partial X^v)_w))
$$
in the case where the dominant weight $\lambda$ is of classical type 
(as introduced in \cite{LS1}, cf. the next definition). Along the way,
we shall retrieve the results of loc. cit., using our basis
$\{p_\pi\}$. In particular, we shall give a geometric characterization
of ``admissible pairs" of loc. cit. (cf. Definition \ref{admissible}
below).

\begin{defn}\label{classical}
Let $\lambda$ be a dominant weight. We say, $\lambda$ is of
classical type if $\langle\lambda, \beta^\vee\rangle\le 2$, for all 
$\beta\in R^+$.
\end{defn}

\medskip

\noindent
{\bf Remarks.}
\begin{enumerate}

\item 
Any dominant weight of classical type is either fundamental, or
a sum of two minuscule fundamental weights.

\item 
$G$ is classical if and only if all fundamental weights of $G$
are of classical type.

\end{enumerate}

\medskip

For the rest of this section, we fix a dominant weight $\lambda$
of classical type.

\begin{proposition}\label{1diml}
Let $v,w \in W^\lambda, v\le w$. Then the $T$--module
$H_w^v(\lambda)$ is at most one--dimensional; further, if non--zero,
then it has the weight $-\frac{1}{2}(w(\lambda)+v(\lambda))$.

As a consequence, the weights of the $T$--module
$H^0(X_w,L_\lambda(-\partial X_w))$ are among the 
$-\frac{1}{2}(w(\lambda)+x(\lambda))$ where $x\leq w$, and the
corresponding weight spaces are one--dimensional.
\end{proposition}

\begin{proof}
Let $p\in H_w^v(\lambda)$. Then $p^2$ belongs to 
$H^0(X_w^v,L_{2\lambda})$, and vanishes of order  
$\ge 2$ along each component of the whole boundary
$(\partial X_w)^v\cup (\partial X^v)_w$. 
On the other hand, the product $p_w p_v$ also belongs to
$H^0(X_w^v, L_{2\lambda})$ and satisfies by Chevalley's formula:
$$
\di(p_w p_v)=\sum_\beta \langle\lambda,\beta^\vee\rangle X^v_{ws_\beta}
+\sum_\gamma \langle\lambda,\gamma^\vee\rangle X^{vs_\gamma}_{w},
$$ 
where $X_{ws_\beta}$ (resp. $X^{vs_\gamma}$) runs over all the
components $X_x$ (resp. $X^y$) of $\partial X_w$ (resp. 
$\partial X^v$) such that $x\ge v$ (resp. $y\le w$). Hence, $p_wp_v$
vanishes of order at most $2$ along each component of 
$(\partial X_w)^v\cup (\partial X^v)_w$ (since $\lambda$ is of
classical type), and nowhere else. Thus, $\frac{p^2}{p_wp_v}$ (a
rational function on $X_w^v$) has no poles. It follows that
$p^2=cp_wp_v, c\in k$, and hence that $p$ 
is unique up to scalars; further, $p$ is either zero or has weight
$\frac{1}{2}$(weight $p_w$+weight $p_v$)
$=-\frac{1}{2}(w(\lambda)+v(\lambda))$.
\end{proof}

As a corollary to the proof of the above Proposition, we have

\begin{lemma}\label{double}
Let $v,w \in W^\lambda, v\le w$. Further, let $H_w^v(\lambda)$ be
non--zero. Then for each divisor $X_{ws_\beta}$ (resp. $X^{vs_\gamma}$)
of $X_w$ (resp. $X^v$) such that ${ws_\beta}\ge v$ (resp. 
${vs_\gamma}\le w$), $\beta$ (resp. $\gamma$) being in $R^+$, we have, 
$\langle\lambda,\beta^\vee\rangle$ 
(resp.~$\langle\lambda,\gamma^\vee\rangle)=2$.
\end{lemma}

We shall denote by $p_{w,v}$ the unique $p_{w,v}^\xi$, if
non--zero (then $p_{w,w}=p_w$). By Proposition \ref{1diml}, 
$p_{w,v}$ lifts to a unique $T$--eigenvector in 
$H^0(X_w,L_\lambda(-\partial X_w))$; we still denote that lift by
$p_{w,v}$. The non--zero $p_{w,v}$, where $v\le w$, form a basis of
$H^0(X_w,L_\lambda(-\partial X_w))$.

Notice that $\frac{p_{w,v}^2}{p_w}$ is a rational section of
$L_\lambda$ on $X_w$, eigenvector of $T$ with weight $-v(\lambda)$,
and without poles by the argument of Proposition \ref{1diml}. This
implies

\begin{lemma}\label{quadratic}
With notations as above, we have $p_{w,v}^2=p_w p_v$ on $X_w$, up to
a non--zero scalar.
\end{lemma}

We now aim at characterizing those pairs $(v,w)$ such that 
$p_{w,v}\neq 0$. For this, we recall some definitions and Lemmas from
\cite{LS1}.

\begin{defn}
Let $X_v$ be a Schubert divisor in $X_w$; further, let $v=s_\alpha w$
where $\alpha\in R^+$. If $\alpha$ is simple, then we say, $X_v$ is a
moving divisor in $X_w$, moved by $\alpha$.
\end{defn}

\begin{lemma}\label{moving} (\cite{LS1} Lemma 1.5.)
Let $X_v$ be a moving divisor in $X_w$, moved by $\alpha$. Let
$X_u$ be any Schubert subvariety of $X_w$. Then either,
$X_u\subseteq X_v$ or $X_{s_\alpha u}\subseteq X_v$.
\end{lemma}

\begin{defn}
Let $v,w \in W^\lambda, v\le w, \ell(v)=\ell(w)-1$; further let
$v=ws_\beta=s_\gamma w$, for some positive roots $\beta,\gamma$. We
denote the positive integer 
$\langle \lambda,\beta^\vee\rangle 
(= \langle v(\lambda),\gamma^\vee\rangle =
-\langle w(\lambda),\gamma^\vee\rangle)$
by $m_\lambda(v,w)$, and refer to it as the Chevalley multiplicity
of $X_v$ in $X_w$ (see \cite{Ch}).
\end{defn}

\begin{lemma}\label{mult} (\cite{LS1} Lemma 2.5.)
Let $v,w \in W^\lambda$ such that $X_v$ is a moving divisor in $X_w$,
moved by $\alpha$. Let $X_u$ be another Schubert divisor in $X_w$. Then
$X_{s_\alpha u}$ is a divisor in $X_v$, and 
$m_\lambda(s_\alpha u,v)=m_\lambda(u,w)$.
\end{lemma}


\begin{defn}
Let $v,w \in W^\lambda$ such that $X_v$ is a divisor in $X_w$.
If $m_\lambda(v,w)=2$, then we shall refer to $X_v$ as a double
divisor in $X_w$.
\end{defn}

By Lemma \ref{double}, if $p_{w,v}\neq 0$, then all Schubert
divisors in $X_w$ that meet $X^v$ are double divisors.

\begin{lemma}\label{dmd} (\cite{LS1} Lemma 2.6.)
Let $u,w \in W^\lambda$ such that $X_u$ is a double divisor in $X_w$. 
Then $X_u$ is a moving divisor in $X_w$.
\end{lemma}


\medskip

\noindent
{\bf{A geometric characterization of Admissible pairs:}} 
Recall (cf.\cite{LS1}):

\begin{defn}\label{admissible}
A pair $(v,w)$ in $W^\lambda$ is called admissible if either
$v=w$ (in which case, it is called a trivial admissible pair), or
there exists a sequence $w=w_1>w_2>\cdots > w_r=v$, such that
$X_{w_{i+1}}$ is a double divisor in $X_{w_{i}}$, i.e.,
$m_\lambda(w_{i+1},w_i)=2$. We shall refer to such a chain as a
double chain.
\end{defn}

We shall give a geometric characterization of admissible pairs (cf.
Proposition \ref{geo} below). First we prove some preparatory
Lemmas.

\begin{lemma}\label{induction} Let $p_{w,v}\neq 0$, then
\begin{enumerate}

\item For any double divisor $X_{s_\alpha w}$ in $X_w$ meeting $X^v$,
we have
$$
p_{s_\alpha w,v}=e_{-\alpha} p_{w,v} \text{ and }
e_{-\alpha}^2 p_{w,v}=0,
$$ 
where $e_{-\alpha}$ is a generator of the Lie algebra of
$U_{-\alpha}$. Further, $p_{s_\alpha w,v}\neq 0$.

\item Likewise, for any double divisor $X^{s_\alpha v}$ in $X^v$
meeting $X_w$, we have
$$
p_{w,s_\alpha v}=e_\alpha p_{w,v} \text{ and }
e_\alpha^2 p_{w,v}=0,
$$ 
where $e_\alpha$ is a generator of the Lie algebra of
$U_\alpha$. Further, $p_{w,s_\alpha v}\neq 0$.

\item The pair $(v,w)$ is admissible.
\end{enumerate}
\end{lemma}

\begin{proof}
(1) Consider the $T$--module $H^0(X_w^v,L_\lambda(-(\partial X^v)_w))$. 
By Proposition \ref{1diml}, it has a basis 
$\{p_{x,v}\vert\,v\leq x\leq w\}$ with corresponding weights
$-\frac{1}{2}(x(\lambda)+v(\lambda))$. Notice that $X_w$ is invariant
under $U_{-\alpha}$ (since $s_\alpha w<w$); hence $X_w^v$ and
$(\partial X^v)_w$ are also $U_{-\alpha}$--invariant. Thus,
$U_{-\alpha}$ acts on $H^0(X_w^v,L_\lambda(-(\partial X^v)_w))$,
compatibly with the $T$--action. The $U_{-\alpha}$--submodule $M$
generated by $p_{w,v}$ is $T$--invariant, with weights of the form
$-\frac{1}{2}(w(\lambda)+v(\lambda))-m\alpha$
for some non--negative integers $m$. But if
$x(\lambda)=w(\lambda)+2m\alpha$, then either $x=w$ and $m=0$, or 
$x=s_\alpha w$ and $m=1$ (by Lemma \ref{double}). 
Hence $M$ is either spanned by $p_{w,v}$, or by $p_{w,v}$ and
$p_{s_\alpha w,v}$. Further, $e_{-\alpha}^2 p_{w,v}=0$.

To complete the proof, it suffices to show that $U_{-\alpha}$ does not
fix $p_{w,v}$. Otherwise, the zero locus of $p_{w,v}$ in $X_w^v$ is
$U_{-\alpha}$--invariant, and hence so is $(\partial X_w)^v$. Thus,
$$
\overline{U_{-\alpha}e_{s_{\alpha w}}}\subseteq (\partial X_w)^v.
$$
But $e_w\in \overline{U_{-\alpha}e_{s_{\alpha w}}}$ (since 
$s_\alpha w<w$) and $e_w\notin (\partial X_w)^v$, a contradiction.

(2) is checked similarly. And (3) follows from (1) together with Lemma
\ref{double}, by induction on $\ell(w)$.
\end{proof}

\begin{lemma}\label{inductionbis}
Let $(v,w)$ be an admissible pair, then $p_{w,v}\neq 0$.
\end{lemma}

\begin{proof}
We argue by induction on $\ell(w)$. We may chose a simple root
$\alpha$ such that $w>s_\alpha w \geq v$ and that $X_{ws_\alpha}$ is a
double divisor in $X_w$. Then 
$\langle w(\lambda),\check\alpha\rangle =-2$, 
and also $p_{s_\alpha w,v}\neq 0$ by the induction hypothesis. The
weight of this vector is 
$$
-\frac{1}{2}(s_\alpha w(\lambda)+v(\lambda))=
-\frac{1}{2}(w(\lambda)+v(\lambda))-\alpha.
$$ 
The scalar product of this weight with $\check\alpha$ being integral,
$\langle v(\lambda),\check\alpha \rangle$ is an even integer. 
Since $\lambda$ is of classical type, it follows that 
$$
\langle v(\lambda),\check\alpha \rangle \in\{2,0,-2\}.
$$
We now distinguish the following three cases:

\medskip

\noindent
{\bf Case 1:} $(v(\lambda),\alpha^\vee)=2$. Then 
$w\geq s_\alpha w, s_\alpha v>v$. 
As a first step, we find a relation between
$H^0(X_w^v,L_\lambda(-(\partial X_w)^v))$ and
$H^0(X_{s_\alpha w}^v,L_\lambda(-(\partial X_{s_\alpha w})^v))$.

Let $G_\alpha$ be the subgroup of $G$ generated
by $U_\alpha$, $U_{-\alpha}$ and $T$; let $B_\alpha=G_\alpha\cap B$. 
Then the derived subgroup of $G_\alpha$ is isomorphic to $SL(2)$ or to
$PSL(2)$, and $G_\alpha/B_\alpha$ is isomorphic to the projective line
$\mP^1$. For a $B_\alpha$--module $M$, we shall denote the associated
$G_\alpha$--linearized locally free sheaf on $G_\alpha/B_\alpha$ by
${\underline{M}}$.

Notice that $X_w$, $X^v$ and hence $X_w^v$ are invariant under
$G_\alpha$, and $(\partial X_w)^v$ is invariant under $B_\alpha$; we
have 
$$
(\partial X_w)^v =
X_{s_\alpha w}^v\cup G_\alpha (\partial X_{s_\alpha w})^v.
$$
Consider the fiber product 
$G_\alpha\times^{B_\alpha} X_{s_\alpha w}^v$ with projection
$$
p:G_\alpha\times^{B_\alpha} X_{s_\alpha w}^v \longrightarrow 
G_\alpha/B_\alpha\simeq\mP^1
$$
and ``multiplication'' map
$$
\psi:G_\alpha\times^{B_\alpha} X_{s_\alpha w}^v \longrightarrow X_w^v.
$$
Then $\psi$ is birational (since it is an isomorphism at $e_w$). 
Further, we have 
$$
(\partial X_w)^v=\psi(X_{s_\alpha w}^v\cup 
G_\alpha\times^{B_\alpha} (\partial X_{s_\alpha w})^v)
$$
where $X_{s_\alpha w}^v$ is the fiber of $p$ at $B_\alpha/B_\alpha$.
By the projection formula, it follows that
$$
L_\lambda(-(\partial X_w)^v) =
\psi_*\psi^*L_\lambda(-X_{s_\alpha w}^v -
G_\alpha\times^{B_\alpha} (\partial X_{s_\alpha w})^v).
$$
This yields an isomorphism
$$
H^0(X_w^v,L_\lambda(-(\partial X_w)^v))\cong 
H^0(G_\alpha/B_\alpha,p_*\psi^*L_\lambda(-X_{s_\alpha w}^v -
G_\alpha\times^{B_\alpha} (\partial X_{s_\alpha w})^v)).
$$
Further, we may identify the $G_\alpha$--linearized sheaf
$p_*\psi^*L_\lambda(-G_\alpha\times^{B_\alpha} 
(\partial X_{s_\alpha w})^v))$ on $G_\alpha/B_\alpha$, to the sheaf 
$\underline{H^0(X_{s_\alpha w}^v,
L_\lambda(-(\partial X_{s_\alpha w})^v))}.$
Therefore, we obtain an exact sequence of $B_\alpha$--modules
$$\displaylines{
0\longrightarrow H^0(X_w^v,L_\lambda(-(\partial X_w)^v))
\longrightarrow 
H^0(G_\alpha/B_\alpha,\underline{H^0(X_{s_\alpha w}^v,
L_\lambda(-(\partial X_{s_\alpha w})^v)})
\hfill\cr\hfill
\longrightarrow
H^0(X_{s_\alpha w}^v,L_\lambda(-(\partial X_{s_\alpha w})^v))
\longrightarrow 0,
\cr}$$
where the map on the right is the ``evaluation'' map (its surjectivity 
follows e.g. from Corollary \ref{surj}.)

Next we analyse the $B_\alpha$--module 
$H^0(X_{s_\alpha w}^v,L_\lambda(-(\partial X_{s_\alpha w})^v))$.
By Proposition \ref{1diml}, its weights have multiplicity one; they are
among the $-\frac{1}{2}(s_\alpha w(\lambda)+x(\lambda))$, where 
$v\leq x\leq s_\alpha w$, and the weight 
$-\frac{1}{2}(w(\lambda)+v(\lambda))-\alpha$
occurs, since $p_{s_\alpha w,v}\neq 0$; its $\alpha$--weight (the
scalar product with $\check\alpha$) is $-2$. 

If $s_\alpha v\nleq s_\alpha w$, then the span $M$ of 
$p_{s_\alpha w,v}$ is invariant under $B_\alpha$. Thus, the 
$T$--module $H^0(G_\alpha/B_\alpha,\underline{M})$ has weights
$-\frac{1}{2}(w(\lambda)+v(\lambda))+(m-1)\alpha$, $m=0,1,2$,
each of them having multiplicitly one. Further, the kernel of the
evaluation map $H^0(G_\alpha/B_\alpha,\underline{M})\longrightarrow M$ 
contains an element of weight $-\frac{1}{2}(w(\lambda)+v(\lambda))$.
By the exact sequence above, this weight occurs in 
$H^0(X_w^v,L_\lambda(-(\partial X_w)^v))$; using Proposition \ref{1diml}
again, it follows that $p_{w,v}\neq 0$. 

On the other hand, if $s_\alpha v\leq s_\alpha w$, then
Lemma \ref{induction} (2) applied to $(s_\alpha w,v)$ yields
$$
p_{s_\alpha w,s_\alpha v}= e_\alpha p_{s_\alpha w,v}\neq 0.
$$
Hence the span $M$ of $p_{s_\alpha w,v}$ and 
$p_{s_\alpha w,s_\alpha v}$ is a non--trivial $B_\alpha$-module with
$\alpha$-weights $-2$ and $0$ (note that 
$e_\alpha p_{s_\alpha w,s_\alpha v}=0$ by weight considerations). 
Thus, we have an isomorphism of $B_\alpha$--modules 
$$
M\cong M_1\otimes M_2,
$$ 
where $M_1$ is a one--dimensional $B_\alpha$--module with
$\alpha$--weight $-1$, and $M_2$ is the standard two--dimensional
$G_\alpha$--module. It follows that the weights of the $T$--module 
$$
H^0(G_\alpha/B_\alpha,\underline{M})\cong 
H^0(G_\alpha/B_\alpha,\underline{M_1})\otimes M_2
$$ 
are exactly $-\frac{1}{2}(w(\lambda)+v(\lambda))-\alpha$, 
$-\frac{1}{2}(w(\lambda)+v(\lambda))+\alpha$
(both of multiplicity one) and
$-\frac{1}{2}(w(\lambda)+v(\lambda))$ (of multiplicity two).
Thus, the kernel of the evaluation map 
$H^0(G_\alpha/B_\alpha,\underline{M})\longrightarrow M$
contains an element of weight $-\frac{1}{2}(w(\lambda)+v(\lambda))$,
and we conclude as above.

\medskip

\noindent
{\bf Case 2:} $\langle v(\lambda),\check\alpha\rangle =0$. Then 
$w>s_\alpha w\geq v=s_\alpha v$, so that $X_w$, $X^v$ and $X_w^v$ are
again invariant under $G_\alpha$, whereas $(\partial X_w)^v$ is
$B_\alpha$--invariant. Arguing as in Case 1, we obtain the same
relation between $H^0(X_w^v,L_\lambda(-(\partial X_w)^v))$ 
and $H^0(X_{s_\alpha w}^v,L_\lambda(-(\partial X_{s_\alpha w})^v))$;
but now the latter $B_\alpha$--module contains the span $M$ of
$p_{s_\alpha w,v}$, as a $B_\alpha$-submodule of $\alpha$--weight
$-1$. As in Case 1, it follows that $p_{w,v}\neq 0$.

\medskip

\noindent
{\bf Case 3:} $\langle v(\lambda),\check\alpha\rangle =-2$. Then
$w>s_\alpha w\geq v>s_\alpha v$, and $X^v$ is a double divisor in
$X^{s_\alpha v}$. Therefore, the pair $(s_\alpha v,s_\alpha w)$ is
admissible. By the induction hypothesis, we have, 
$p_{s_\alpha v,s_\alpha w}\neq 0$. Then Case 1 applies to the pair
$(s_\alpha v,w)$ and yields $p_{w,s_\alpha v}\neq 0$. Further, $X^v$
is a double divisor in $X^{s_\alpha v}$. Hence by Lemma
\ref{induction} (2) applied to $(w,s_\alpha v)$, we obtain  
$p_{w,v}\neq 0$. 
\end{proof}

Now combining Lemmas \ref{double}, \ref{induction} and
\ref{inductionbis}, we obtain 

\begin{proposition}\label{geo}
Let $v,w \in W^\lambda, v\le w$. Then the pair $(v,w)$ is admissible
if and only if $p_{w,v}$ is non--zero. In this case, every chain from
$v$ to $w$ is a double chain.
\end{proposition}

\section{Standard monomials for sums of weights of classical type}

In this section, we obtain a standard monomial basis for
$H^0(X_w^v,L_{\lambda_1+\cdots+\lambda_m})$, where $X_w^v$ is a
Richardson variety in $G/P$, and $\lambda_1,\dots,\lambda_m$ are
dominant characters of classical type of $P$ (in the sense of
Definition \ref{classical}). 

We begin with the case where $m=1$; we shall need a definition, and a
result of Deodhar (\cite{LS1} Lemmas 4.4 and 4.4') on the Bruhat
ordering.

\begin{defn}\label{lambda}
Let $w\in W^P$ and let $\lambda$ be a dominant character of
$P$. We say that $x\in W^P$ is $\lambda$--maximal in $w$
(resp.~$\lambda$--minimal on $w$) if $xy\le x$ for any $y\in W_\lambda$
such that $xy\in W^P$ and $xy\le w$ (resp.~if $xy\ge x$ for any 
$y\in W_\lambda$ such that $xy\in W^P$ and $xy\ge w$).
\end{defn}

\begin{lemma}\label{deodhar}
Let $w\in W$ and $x\in W^\lambda$ such that 
$x\le w(\lambda)$ (resp.~$x\geq v(\lambda)$). Then
the set $\{y\in W_\lambda~\vert~ xy\leq w\}$
(resp.~$\{y\in W_\lambda~\vert~ w\leq xy\}$ 
admits a unique maximal (resp.~minimal) element.
\end{lemma}

We shall also need the following consequences of this result.

\begin{lemma}\label{maximal}

\begin{enumerate}

\item
Let $w\in W^P$ and $x\in W^\lambda$ such that $x\le w(\lambda)$
(resp.~$x\ge w(\lambda)$). Then $x\in W/W_\lambda$ admits a unique
lift $\tilde x\in W^P$ such that $\tilde x$ is $\lambda$--maximal in $w$
(resp.~$\lambda$--minimal on $w$).

\item
Let $v\le w\in W^P$, then $v$ is $\lambda$--maximal in $w$ 
(resp.~$w$ is $\lambda$--minimal in $v$) if and only if 
$(\partial_\lambda X^v)_w=(\partial X^v)_w$ 
(resp. $(\partial_\lambda X_w)^v = (\partial X_w)^v$).

\end{enumerate}

\end{lemma}

\begin{proof}
(1) Let $x\le w(\lambda)$. By Lemma \ref{deodhar}, the set 
$\{y\in W_\lambda ~\vert~ xy\le w\}$
admits a unique maximal element that we still denote by $y$. Let 
$\tilde x$ be the representative in $W^P$ of $xy\in W$, then 
$\tilde x(\lambda)=xy(\lambda)= x(\lambda)$. Further, if we have 
$\tilde x z\leq w$ for some $z\in W_\lambda$ 
such that $\tilde x z\in W^P$, then we can write $\tilde x z = x u$ where 
$u\in W_\lambda$. Since $xu\le w$, we have $u\le y$ and hence 
$xu\le xy$  (since $x\in W^\lambda$ and $u,y\in W_\lambda$). But
$xu=\tilde x z\in W^P$, so that $\tilde x z\le \tilde x$. This proves the
assertion concerning $\lambda$--maximal elements, and hence the dual
assertion concerning $\lambda$--minimal elements.

(2) If $(\partial_\lambda X^v)_w \neq (\partial X^v)_w$, then there
exists $y\in W_\lambda$ such that $v<  vy \le w$ and
$\ell(vy)=\ell(v)+1$. Thus, $v$ is not $\lambda$--maximal in $w$.

Conversely, if $v$ is not $\lambda$--maximal in $w$, then 
$v < \tilde v \le w$ where $\tilde v\in vW_\lambda$ is
$\lambda$--maximal in $w$. Hence there exists $y\in W_\lambda$ such
that $v< vy \le \tilde v \le w$ and $\ell(vy)=\ell(v)+1$. Now
$X_w^{vy}$ is contained in $(\partial X^v)_w$ but not in 
$(\partial_\lambda X^v)_w$.
\end{proof}

Now we consider the $T$--module $H^0(X_w^v,L_\lambda)$, where
$v\le w\in W^P$ and $\lambda$ is a dominant character of $P$, not
necessarily $P$--regular. Notice that the diagram
$$\CD
H^0(X_{w(\lambda)},L_\lambda)@>>>
H^0(X_{w(\lambda)}^{v(\lambda)},L_\lambda)\\
@VVV @VVV\\
H^0(X_w,L_\lambda)@>>> H^0(X_w^v,L_\lambda)
\endCD$$
is commutative, where the horizontal (resp.~vertical) maps are
restrictions (resp.~pull--backs). Further, both restrictions are
surjective by Proposition \ref{vanishing1}; and the pull--back on the
left is an isomorphism, since the natural map
$f:X_w\to X_{w(\lambda)}$ satisfies
$f_*\cO_{X_w}=\cO_{X_{w(\lambda)}}$. Thus, we may regard the
$T$--module 
$H^0(X_w^v,L_\lambda)$ as a quotient of
$H^0(X_{w(\lambda)}^{v(\lambda)},L_\lambda)$. 

Likewise, by using the commutative diagram 
$$\CD
H^0(X_{w(\lambda)},L_\lambda(-\partial X_{w(\lambda)}))@>>>
H^0(X_{w(\lambda)}^{v(\lambda)},
L_\lambda(-(\partial X_{w(\lambda)})^{v(\lambda)}))\\
@VVV @VVV\\
H^0(X_w,L_\lambda(-\partial_\lambda X_w))@>>> 
H^0(X_w^v,L_\lambda(-(\partial_\lambda X_w)^v))
\endCD$$
and Corollary \ref{surj}, we may regard the $T$--module
$H^0(X_w^v,L_\lambda(-(\partial_\lambda X_w)^v))$ as a quotient of 
$H^0(X_{w(\lambda)}^{v(\lambda)},
L_\lambda(-(\partial X_{w(\lambda)})^{v(\lambda)}))$. The latter has
been described in Section 6, in the case that $\lambda$ is of
classical type: it has a basis consisting of the
$p_{w(\lambda),x(\lambda)}$ where $x(\lambda)\in W^\lambda$,
$v(\lambda)\le x(\lambda)\le w(\lambda)$ and the pair 
$(x(\lambda),w(\lambda))$ is admissible. From that description we
shall deduce

\begin{proposition}\label{oneweight}
Let $v\le w\in W^P$ and let $\lambda$ be a dominant character of
classical type of $P$.

\begin{enumerate}

\item
The space
$H^0(X_w^v,L_\lambda(-(\partial_\lambda X_w)^v-(\partial X^v)_w))$
is spanned by $p_{w(\lambda),v(\lambda)}$, if $v$ is $\lambda$--maximal
in $w$; otherwise, this space is zero.

\item
The $p_{w(\lambda),x(\lambda)}$ where $x\in W^P$ and 
$v\le x\le w$, form a basis of the space
$H^0(X_w^v,L_\lambda(-(\partial_\lambda X_w)^v))$.

\item
The $p_{w(\lambda),x(\lambda)}$ where $x\in W^P$ is $\lambda$--minimal
on $v$, and $w$ is $\lambda$--minimal on $x$, form a basis of
$H^0(X_w^v,L_\lambda(-(\partial X_w)^v))$.

\item
The $p_{y(\lambda),x(\lambda)}$ where $x,y\in W^P$ and 
$v\le x\le y\le w$, form a basis of $H^0(X_w^v,L_\lambda)$.

\end{enumerate}
\end{proposition}

\begin{proof}
(1) Assume that 
$H^0(X_w^v,L_\lambda(-(\partial_\lambda X_w)^v-(\partial X^v)_w))$
contains a non--zero element $p$. Then, by the argument of Proposition 
\ref{1diml}, $\frac{p^2}{p_{w(\lambda)}p_{v(\lambda)}}$ is a rational
function on $X_w^v$, without poles; further, it vanishes identically
on $(\partial X^v)_w - (\partial_\lambda X^v)_w$, since the zero locus
of $p_{v(\lambda)}$ is $(\partial_\lambda X^v)_w$. It follows that
$p^2$ is a constant multiple of $p_{w(\lambda)}p_{v(\lambda)}$, and
that $(\partial X^v)_w  = (\partial_\lambda X^v)_w$. Hence $p$ is a
constant multiple of $p_{w(\lambda),v(\lambda)}$, and $v$ is
$\lambda$--maximal in $w$ (by Lemma \ref{maximal}).

Conversely, let $v$ be $\lambda$--maximal in $w$; then 
$(\partial_\lambda X^v)_w=(\partial X^v)_w$. Thus,
$p_{w(\lambda),v(\lambda)}$ vanishes identically on 
$(\partial_\lambda X_w)^v\cup(\partial X^v)_w$. Further,
$p_{w(\lambda),v(\lambda)}\ne 0$ on $X_w^v$, since 
$p_{w(\lambda),v(\lambda)}^2=p_{w(\lambda)}p_{v(\lambda)}$ on $X_w$
(by Lemma \ref{quadratic}).

(2) By  Proposition \ref{1diml}, the space
$H^0(X_w^v,L_\lambda(-(\partial_\lambda X_w)^v))$ is spanned by the
images of the $p_{w(\lambda),x(\lambda)}$ where 
$v(\lambda)\le x(\lambda)\le w(\lambda)$. Further,
$p_{w(\lambda),x(\lambda)}^2 = p_{w(\lambda)} p_{x(\lambda)}$ on
$X_w$. Using Lemma \ref{extremal}, we see that 
$p_{w(\lambda),x(\lambda)}$ is non--zero on $X_w^v$ if and only if
$x(\lambda)$ has a representative $x\in W^P$ such that $v\le x\le w$.

(3) $H^0(X_w^v,L_\lambda(-(\partial X_w)^v))$
is a $T$--stable subspace of 
$H^0(X_w^v,L_\lambda(-(\partial_\lambda X_w)^v))$; thus, it is spanned
by certain $p_{w(\lambda),x(\lambda)}$ where $v\le x\le w$. By Lemma
\ref{quadratic}, the zero locus $(p_{w(\lambda),x(\lambda)}=0)$ in
$X_w^v$ equals $(p_{x(\lambda)}=0)\cup(\partial_\lambda X_w)^v$. Hence
$p_{w(\lambda),x(\lambda)}$ belongs to 
$H^0(X_w^v,L_\lambda(-(\partial X_w)^v))$ 
if and only if $p_{x(\lambda)}$ vanishes identically on 
$(\partial X_w)^v - (\partial_\lambda X_w)^v$. By Lemma
\ref{extremal}, this amounts to: $x(\lambda)$ admits no lift $x'$ such
that $v\le x'\le wy$ for some $y\in W_\lambda$, $wy<y$,
$\ell(wy)=\ell(w)-1$. Let $\tilde x$ be the lift of $x(\lambda)$ that is
$\lambda$--minimal on $v$, then the preceding condition means that $w$
is $\lambda$--minimal on $\tilde x$.

(4) By Proposition \ref{filtration1}, we obtain a basis of the space
$H^0(X_w^v,L_\lambda)$ by choosing a basis of
$H^0(X_w^x,L_\lambda(-(\partial X^x)_w))$ for each $x\in W^P$ such
that $v\le x\le w$, and lifting this basis to $H^0(X_w^v,L_\lambda)$
under the (surjective) restriction map
$H^0(X_w^v,L_\lambda)\to H^0(X_w^x,L_\lambda)$. Together with (3), it
follows that a basis of $H^0(X_w^v,L_\lambda)$ consists of the
$p_{y(\lambda),x(\lambda)}$, where $y$ is $\lambda$--maximal in $w$,
$x$ is $\lambda$--maximal in $y$, and $v\le x$. But given any 
$x',y'\in W^P$ such that $v\le x'\le y'\le w$, we have
$v\le x\le y\le w$ and
$p_{y'(\lambda),x'(\lambda)}=p_{y(\lambda),x(\lambda)}$, where $x$
(resp.~$y$) is the representative of $x(\lambda)$ that is
$\lambda$--maximal in $w$ (resp.~$x$).
\end{proof}

\begin{defn}
Let $\lambda_1,\dots,\lambda_m$ be dominant characters of classical
type of $P$. Let $\pi_i=(w_i,v_i)$ where  
$v_i\le w_i\in W^{\lambda_i}$ for $1\le i\le m$. Then the sequence
$\underline{\pi}=(\pi_1,\dots,\pi_m)$ is standard if there exist lifts 
$\tilde w_i,\tilde v_i$ in $W^P$ for $1\le i\le m$, such that
$$
\tilde v_m\le \tilde w_m\le \cdots \le \tilde v_1\le \tilde w_1.
$$ 
The monomial 
$$
p_{\underline{\pi}}=p_{w_1(\lambda_1),v_1(\lambda_1)}\cdots
p_{w_m(\lambda_m),v_m(\lambda_m)}
\in H^0(G/P,L_{\lambda_1+\cdots+\lambda_m})
$$ 
is called standard as well.

Further, let $v,w\in W^P$ such that $v\le w$; then $\underline{\pi}$
is standard on $X_w^v$ if there exist lifts as above, such that
$$
v\le \tilde v_m\le \tilde w_m\le \cdots
\le \tilde v_1\le \tilde w_1\le w.
$$
The restriction of $p_{\underline{\pi}}$ to $X_w^v$ is called a
standard monomial on $X_w^v$; it is a 
$T$--eigenvector in $H^0(X_w^v,L_{\lambda_1+\cdots+\lambda_m})$. 
\end{defn}

Notice that there is no loss of generality in assuming that $\tilde v_m$
is $\lambda_m$--minimal on $v$, and that $\tilde w_m$ is 
$\lambda_m$--minimal on $\tilde v_m$. 

Now the argument of Proposition \ref{products}, together with
Proposition \ref{oneweight} and induction on $m$, yields the following
partial generalization of Corollary \ref{monomials}.

\begin{theorem}\label{mixed}
Let $v\le w\in W^P$ and let $\lambda_1,\dots,\lambda_m$ be dominant
characters of $P$. If $\lambda_1,\dots,\lambda_m$ are of classical
type, then the standard monomials on $X_w^v$ form a basis for
$H^0(X_w^v,L_{\lambda_1+\cdots+\lambda_m})$. 
\end{theorem}

\noindent
{\bf Remarks.} 

\begin{enumerate}

\item
In particular, Theorem \ref{mixed} applies to $P=B$ if all fundamental
weights are of classical type, that is, if $G$ is classical. Thereby,
we retrieve all results of \cite{LS1}.

\item 
The second assertion of Corollary \ref{monomials} does
not generalize to this setting, that is, there are examples of
standard monomials on $G/P$ which are not standard on $X_w^v$, but
which restrict non--trivially to that subvariety. 

Specifically, let $G=SL(3)$ with simple reflections $s_1,s_2$ and
fundamental weights $\omega_1,\omega_2$. Then one may check that the
monomial
$$
p_{s_1(\omega_1)} p_{s_2(\omega_2)}\in H^0(G/B,L_{\omega_1+\omega_2})
$$ 
is standard on $G/B$ and restricts non--trivially to $X_{s_2s_1}$, but
is not standard there.
\end{enumerate}


\begin{thebibliography}{100}

\bibitem{B} {\sc M.~Brion}: Positivity in the Grothendieck group of
complex flag varieties, preprint available at math.AG/0105254.

\bibitem{BP} {\sc M.~Brion} and {\sc P. Polo}: Large Schubert
varieties, {\it Represent. Theory} {\bf 4} (2000), 97--126.

\bibitem{Ch} {\sc C.~Chevalley}: Sur les d\'ecompositions cellulaires
des espaces $G/B$ (with a foreword by A. Borel), 
\textit{Proc. Sympos. Pure Math.} {\bf 56}, Part 1, Algebraic Groups
and their Generalizations: Classical Methods (University Park, PA;
1991), \textit{Amer. Math. Soc.}, Providence, RI (1994), 1--23.

\bibitem{D} {\sc V.~V.~Deodhar}: On some geometric aspects of Bruhat 
orderings. I. A finer decomposition of Bruhat cells. 
{\it Invent. Math.} {\bf 79} (1985), 499--511.

\bibitem{H} {\sc W. V. D. Hodge}:  Some enumerative results in the
theory of forms, {\it Proc. Camb. Phil. Soc.} {\bf 39} (1943), 22-30.

\bibitem{KK} {\sc B.~Kostant} and {\sc S.~Kumar}:
$T$-equivariant $K$-theory of generalized flag varieties, 
{\it J. Differential Geom.} {\bf 32} (1990), 549--603.

\bibitem{Kn} {\sc A.~Knutson}: A Littelmann-type formula for
Duistermaat-Heckman measures, {\it Invent. math.} {\bf 135} (1999),
185--200.

\bibitem{LL} {\sc V.~Lakshmibai} and {\sc P.~Littelmann}:
Richardson varieties and equivariant $K$-theory (work in progress, 2001).

\bibitem{LS1} {\sc V.~Lakshmibai} and {\sc C.S.~Seshadri}: 
Geometry of G/P-V, {\it J. Alg}, {\bf 100} (1986), 462--557.

\bibitem{LS2} {\sc V.~Lakshmibai} and {\sc C.S.~Seshadri}: 
Standard monomial theory, {\it Proc. Hyderabad Conference on Algebraic
Groups},(S. Ramanan et al., eds.), Manoj Prakashan, Madras (1991),
279--323.

\bibitem{Li1}{\sc P.~Littelmann}: A Littlewood-Richardson formula for 
symmetrizable Kac-Moody algebras, {\it Invent. Math.} {\bf 116}
(1994), 329--346. 

\bibitem{Li2} {\sc P.~Littelmann}: Contracting modules and standard 
monomial theory, {\it J. Amer. Math. Soc.} {\bf 11} (1998), 551--567.

\bibitem{MK} {\sc V.~Mehta} and {\sc W.~van der Kallen}: On a
Grauert-Riemenschneider theorem for Frobenius split varieties in
characteristic $p$, {\it Invent. math.} {\bf 108} (1992), 11--13.

\bibitem{RR} {\sc S.~Ramanan} and {\sc A.~Ramanathan}: 
Projective normality of flag varieties and Schubert varieties,
{\it Invent. math.} {\bf 79} (1985), 217--234.

\bibitem{Ram1} {\sc A.~Ramanathan}: Schubert varieties are
arithmetically Cohen-Macaulay, {\it Invent. math.} {\bf 80} (1985), 
283--294.

\bibitem{Ram2} {\sc A.~Ramanathan}: Equations defining Schubert
varieties and Frobenius splitting of diagonals, {\it Pub. Math.
IHES} {\bf 65} (1987), 61--90.

\bibitem{Ric} {\sc R.~W.~Richardson}: Intersections of double cosets
in algebraic groups, {\it Indag. Math. (N. S.)} {\bf 3} (1992),
69--77.

\end{thebibliography}
\end{document}